\documentclass[leqno,twoside, 11pt]{amsart}
\usepackage{latexsym,esint}
\usepackage{graphicx}
\usepackage{color}
\usepackage{verbatim}

\theoremstyle{plain}

\numberwithin{equation}{section} \numberwithin{figure}{section}

\newcommand{\dps}{\displaystyle}

\newtheorem{theorem}{Theorem}[section]
\newtheorem{lemma}[theorem]{Lemma}

\title[On the F\"oppl-von K\'arm\'an theory for prestrained films]{On the F\"oppl-von K\'arm\'an theory for \\elastic prestrained films with varying thickness}
\author{Hui Li}
\address{Changzhou Vocational Institute of Industry Technology, School of Information Engineering,
28 Minxin Middle Rd., Changzhou, Jiangsu, China 213164}
\email{lihui@ciit.edu.cn} 

\begin{document}

\begin{abstract}
    We derive the variational limiting theory of thin films, parallel to the F\"oppl-von K\'arm\'an theory in the nonlinear elasticity, for films that have been prestrained and whose thickness is a general non-constant function. Using $\Gamma$-convergence, we extend the existing results to the variable thickness setting, calculate the associated Euler-Lagrange equations of the limiting energy, and analyze convergence of equilibria. The resulting formulas display the interrelation between deformations of the geometric mid-surface and components of the growth tensor. 
\end{abstract}

\maketitle

\section{Introduction}
The use of the notion of $\Gamma$-convergence in studying elastic thin plates has been first proposed in the mid-1990s \cite {LeDret-Raoult_1996, LeDret-Raoult_1995}, and has rapidly developed in the past thirty years. On the one hand, various $2$ dimensional models have been rigorously derived from the theory of $3$d nonlinear elasticity \cite{Friesecke-James-Muller_2002, Friesecke-James-Muller_2006, Friesecke-James-Mora-Muller_2003, LeDret-Raoult_1996, Lewicka-Li_2015, Lewicka-Mora-Pakzad_2009, Lewicka-Pakzad_2013}, while on the other hand, non-Euclidean elasticity of plates and shells has successfully attempted describing the phenomenon of morphogenesis, with prestrained films as its research objects (see the recent monograph \cite{Lewicka_2023} and references therein).

\begin{figure}[h]
\includegraphics[width = 0.6 \textwidth]{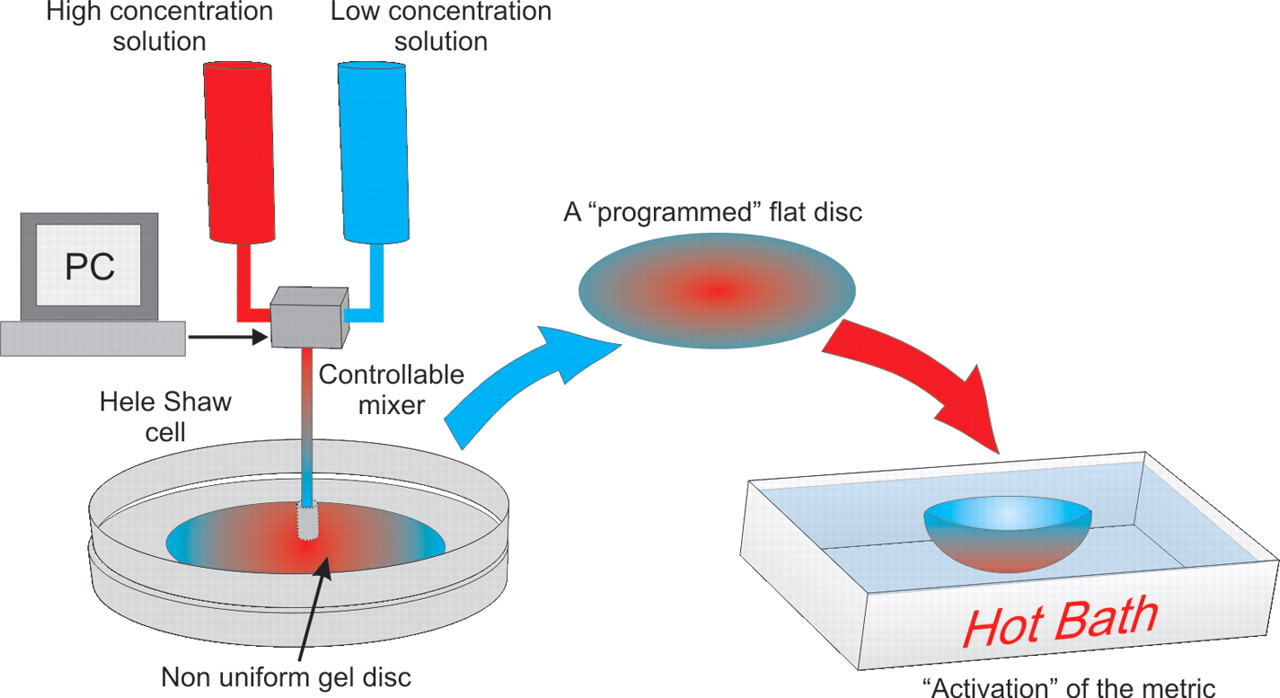}\hspace{6mm}
\includegraphics[width = 0.3 \textwidth]{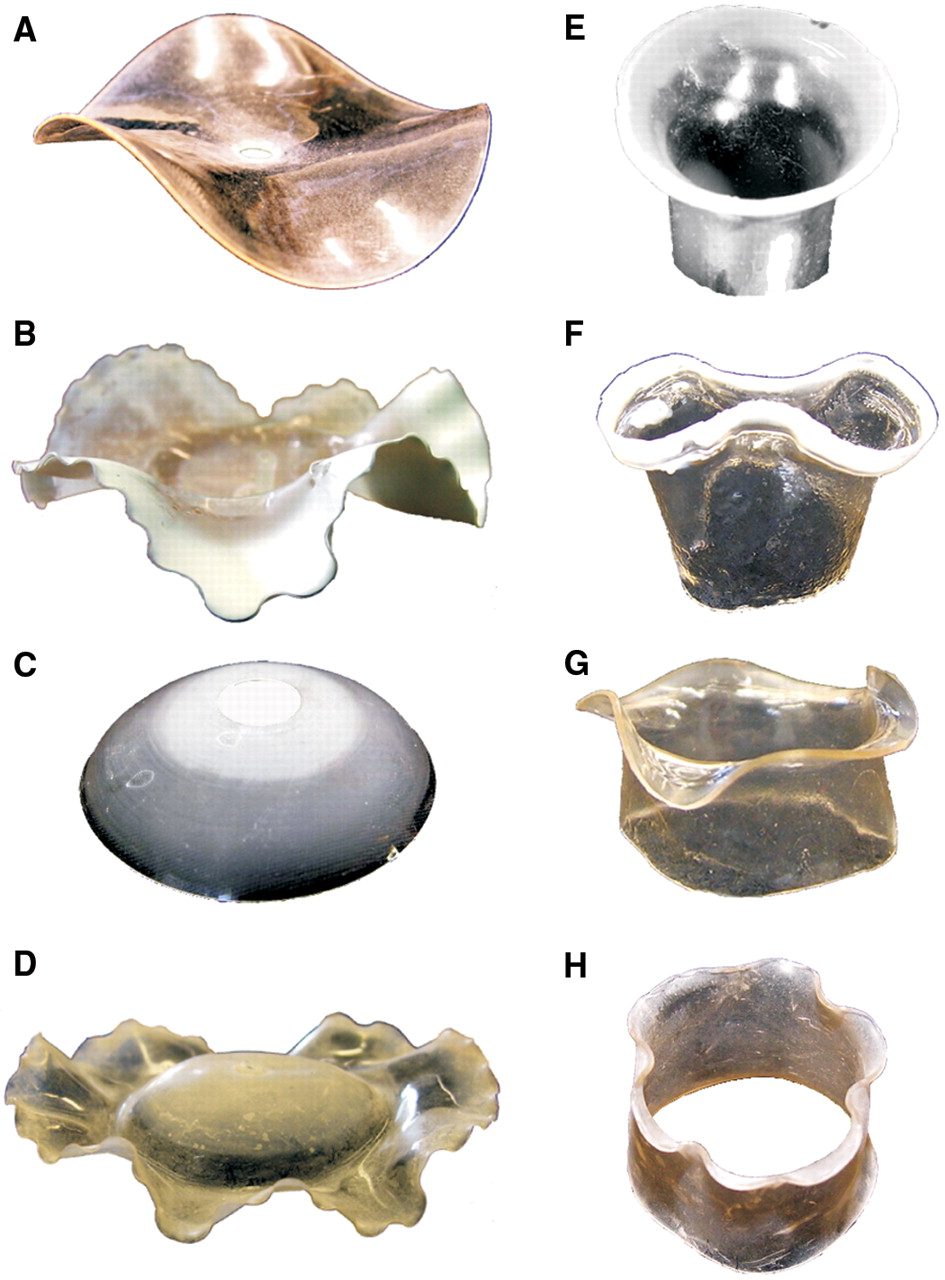}
\caption{Imposing an incompatible target metrics a sheets of NIPA gels. The experiment (on the left) and the obtained film shapes (on the right) in \cite{Klein-Efrati-Sharon_2007}}.
\label{experiment}
\end{figure}

The simple morphogenesis principle, as depicted in Figure \ref{experiment}, proposes that a local heterogeneous incompatibility of strains, represented by an incompatible Riemannian metric $G^h$, posed on a thin referential configuration $S^h$, results in the local elastic stresses \cite{Efrati-Sharon-Kupferman_2009, Klein-Efrati-Sharon_2007}. Thus prestrained films are ubiquitous in nature and engineering applications, such as: growing tissues, plastically strained sheets, swelling or shrinking gels, petals and leaves of flowers, atomically thin graphene layers, to mention a few. In order to fully relieve the tension, $S^h$ strives to realize $G^h$ and settles with a shape, in a sense, closest to the isometric realization of $G^h$. 

The analytical set-up for the non-Euclidean elasticity of thin films is as follows. We assume $S \subset \mathbb{R}^3$ to be a 2d surface, and for each small $h>0$ we pose:
\begin{equation}\label{defSh}
S^h = \left\{z = z' + t \vec n(z')\mid z' \in S, -g_1^h(z') < t < g_2^h(z')\right\},
\end{equation}
where $\vec n$ is the unit normal to $S$ and $g_i^h \sim h$ for $i=1,2$ are scalar positive functions on $S$. Let $G^h$ be a Riemannian metric on $S^h$ and let $u^h \in W^{1, 2}(S^h; \mathbb R^3)$ represent a deformation of $S^h$. We set:
\begin{equation}\label{energygrowthI}
I^h_g(u^h) = \frac{1}{h} \int_{S^h} W\left(\nabla u^h (G^h)^{-1/2}\right)~\mathrm{d}z,
\end{equation}
where $(G^h)^{-1/2}$ is the inverse of the unique symmetric positive definite square root of $G^h$, and $W: \mathbb R^{3 \times 3} \to \mathbb R_+$ is the given energy density function satisfying the following properties of frame indifference, normalization, non-degeneracy, and local regularity:
\begin{equation}\label{w-ass}
\begin{minipage}{12cm}
\begin{itemize}
\item [(i)] $W(RF) = W(F)$, for all $R \in SO(3)$ and $F \in \mathbb R^{3 \times 3}$.
\item [(ii)] $W(\mbox{Id}) = 0$.
\item [(iii)] $W(F) \geq c \; \mbox{dist}^2 (F, SO(3))$ with $c$ being a positive constant.
\item [(iv)] $W$ is $\mathcal C^2$ in a $\delta$-neighborhood of $SO(3)$.
\end{itemize}
\end{minipage}
\end{equation}
We point out that for deformations with gradient close to $SO(3)$, condition (iii) above makes $I^h_g(u^h)$ in (\ref{energygrowthI}) comparable to the functional $\bar I^h(u^h)$ defined as:
\[
\bar I^h(u^h) = \frac{1}{h}\int_{S^h}\mbox{dist}^2\left(\nabla u^h(G^h)^{-1/2}, SO(3)\right)~\mathrm{d}z,
\]
and measuring how well the metric $G^h$ is realized by the deformation $u^h$. Here, $\mbox{dist}(\cdot, SO(3))$ is the distance of a $3\times 3$ matrix from the (compact) special orthogonal group $SO(3)$.

The theory of dimension reduction explores the asymptotic behaviour of the energy functional $I^h_g$ when the thickness parameter $h \to 0$, by first determining the scaling exponent $\beta$ such that $\inf I^h_g \sim h^{\beta}$, then deriving the $\Gamma$-limit $\mathcal I_{\beta}$ of $h^{-\beta}I^h_g$.
We now briefly review the literature corresponding to plates with uniform thickness i.e. $S=\Omega \subset \mathbb R^2$ and $g_1^h=g_2^h = h/2$. The case $\beta \geq 2$ and $G^h=G(z')$ has been discussed in papers \cite{Bhattacharya-Lewicka-Schaffner_2016, Lewicka-Pakzad_2011}. In \cite{Lewicka-Raoult-Ricciotti_2017} it has been shown that if $\beta > 2$, then $\inf I^h \leq h^4$ which further corresponds to the specific condition on the Riemann curvatures 
$\{R_{12,ab}\}_{a,b = 1, 2, 3} = 0$ on $\Omega$. Moreover, if $\beta > 4$, then $\inf I^h_g \leq h^6$, arising when all curvatures satisfy $R(G)=0$ on $\Omega$. The paper \cite{Lewicka_2020} extended these results to having $G^h = G \in \mathcal C^{\infty}(\bar{\Omega}^1, \mathbb R^{3\times 3}_{sym, +})$ variable in the normal direction, and proved that the order of $\inf I^h_g$ relative to  $h$ can only be even, i.e. $\inf I^h_g\sim h^{2n}$, obtaining all $\Gamma$-limits in such infinite hierarchy $\{\mathcal{I}_{2n}\}_{n \geq 1}$ of prestrained thin plates. In comparison, the hierarchy of plate models in classical case nonlinear elsticity presented in \cite{Friesecke-James-Muller_2006}, contains only four such limiting objects: the Kirhchoff, the nonlinear Kirhchoff, the von K\'arm\'an and the linear elasticity. In paper \cite{Lewicka-Lucic_2020} metrics $G^h$ with more the pronounced oscillatory nature are studied, 
while the case of even more general structure of $G^h$ under the assumption of being close to the single immersable metric $\mathrm{Id}_3$, has been discussed in \cite{Lewicka-Mahadevan-Pakzad_2014, Lewicka-Ochoa-Pakzad_2015, Lewicka-Mahadevan-Pakzad_2017}. For non-Euclidean shells, paper \cite{Li_2013} derived the Kirchhoff theory for prestrained shells with the metric  invariant in the normal direction. In the abstract setting of Riemannian manifolds, general results have been also presented \cite{Kuferman-Maor_2014, Kupferman-Soloman-Jake_2014, Maor-Shachar_2019}.
When $\beta < 2$, although no systematic results are available so far, there have been various studies of the emerging patterns in the context of: compression-driven blistering \cite{BenBelgacem-Conti-DeSimone-Muller_2000, BenBelgacem-Conti-DeSimon-Muller_2002, Jin-Sternberg_2001}, buckling \cite{Gemmer-Sharon-Shearman-Venkataramani_2016, Gemmer-Venkataramani_2013, Gemmer-Venkataramani_2011}, origami patterns \cite{Conti-Maggi_2008, Venkataraman_2004}, conical singularities \cite{Muller-Olbermann_2014, Olbermann_2019, Olbermann_2016} and coarsening patterns \cite{Bella-Kohn_2014, Bella-Kohn_2017}. 

All studies mentioned above concern the uniform thickness scenario. However, in both nature and engineering, plates or shells with varying thickness are more common. Although some results exist for the classical nonlinear elasticity \cite{Lewicka-Mora-Pakzad_2009, Li_2017}, little is known in case of the nontrivial prestrain. In the present paper, we will thus address the varying thickness situation as in (\ref{defSh}) for non-Euclidean plates, i.e. $S =   \Omega \subset \mathbb{R}^2$, with $g_1^h, g_2^h$ satisfying: 
\begin{equation}\label{g1hg2hI}
\lim_{h \to 0}\frac{g_1^h}{h} = g_1 \quad \mbox{and} \quad \lim_{h \to 0}\frac{g_2^h}{h} = g_2 \quad \mbox{in}~\mathcal C^1(\bar \Omega),
\end{equation} 
where $g_1, g_2$ are two positive $\mathcal C^1$ functions on $\bar{\Omega}$. To be more experimentally relevant, this paper chooses the growth tensor as in \cite{Lewicka-Mahadevan-Pakzad_2011} and extends the results therein which are the rigorous analytical counterparts of the asymptotic expansion argument in \cite{Liang-Mahadevan_2009}. 

We also derive the resulting Euler-Lagrange equations, generalizing those obtained in \cite{Lewicka-Mahadevan-Pakzad_2011}. 
Finally, under additional physical conditions (\ref{m-w-ass}) for the elastic energy density $W$, we establish convergence of equilibria (rather than only of minimizers) i.e. convergence of critical points of the discussed 3$d$ non-Euclidean energies to the critical points of the corresponding $\Gamma$-limiting energy derived in this paper. Prior studies of such convergence, in case of the classical plates/shells theories appeared in \cite{Muller-Pakzad_2008, Lewicka_2011, Mora-Scardia_2012, Lewicka-Li_2015}, however the prestrained case has not been addressed so far. 

\section{An overview of the main results}\label{Brief_FvK}
We consider a sequence of $3$d thin plates: 
\begin{equation}\label{VShellFvK}
\Omega^h = \{x=(x', x_3)\mid x'\in \Omega, x_3 \in (-g_1^h(x'), g_2^h(x'))\},
\end{equation}
where $\Omega \subset \mathbb R^2$ is an open, bounded, simply connected domain and $g_1^h, g_2^h \in \mathcal{C}^1(\bar \Omega)$ satisfy (\ref{g1hg2hI}). It is convenient to define the universal rescaled domain $\Omega^{\ast}$ in:
\begin{equation}\label{CommonDomainFvK}
\Omega^{\ast} = \{(x', x_3)\mid x' \in \Omega, \; x_3 \in (-1/2, 1/2)\},
\end{equation}
and the change of variable $s^h(x', \cdot):(-1/2, 1/2)\to 
(-g_1^h(x'), g_2^h(x'))$ as:
\begin{equation}\label{ChangeVarFvK}
s^h(x', x_3) = \big(g_1^h(x') + g_2^h(x')\big) x_3 + \frac{1}{2} \big(g_2^h(x')- g_1^h(x')\big).
\end{equation}
Each $\Omega^h$ undergoes an instantaneous 
growth, due to $a^h = [a_{ij}^h]: \Omega^h \to \mathbb R^{3 \times 3}$ of the form:
\begin{equation}\label{GrowthTensorFvK}
a^h(x', x_3) = \mbox{Id}_3 + h^2 \epsilon_g(x') + hx_3 \kappa_g(x'),
\end{equation}
where $\epsilon_g, \kappa_g: \bar \Omega \to \mathbb R^{3 \times 3}$ are two given smooth matrix fields.
For each deformation $u^h \in W^{1, 2}(\Omega^h, \mathbb R^3)$, its elastic energy is now determined similarly to (\ref{energygrowthI}) in:
\begin{equation}\label{VElasticEnergyFvK}
I^h(u^h) = \frac{1}{h} \int_{\Omega^h} W\big(\nabla u^h (a^h)^{-1}\big)\mbox{d}x,
\end{equation}
where the stored energy density $W: \mathbb R^{3 \times 3} \to[0, \infty]$ is as in (\ref{w-ass}). As mentioned in \cite{Lewicka-Mahadevan-Pakzad_2011}, when the energy density $W$ is isotropic, the functional in (\ref{VElasticEnergyFvK}) reduces to (\ref{energygrowthI}) with $G^h = (a^h)^T a^h$.
Combining techniques in \cite{Lewicka-Mahadevan-Pakzad_2011, Li_2017}, in section \ref{VariableT_FvK1}, we derive the limiting energy of $I^h$ as $h\to 0$, as:
\begin{equation}\label{Ig_intro}
\begin{split}
\mathcal I_g(v, w) = & \; \frac{1}{2}\int_{\Omega}(g_1 + g_2)\mathcal Q_2\Big(\mbox{sym} \nabla w + \frac{1}{2} \nabla v \otimes \nabla v - (\mbox{sym} \;\epsilon_g)_{2 \times 2} - \frac{1}{2}(g_2 - g_1)(\mbox{sym}\;\kappa_g)_{2 \times 2}\\
& \qquad\qquad\qquad\qquad  + \frac{1}{2}\mbox{sym}(\nabla v \otimes \nabla(g_2 - g_1))\Big)~\mbox{d}x' \\ & + \frac{1}{24}\int_{\Omega}(g_1 + g_2)^3 \mathcal Q_2 \left(\nabla^2 v + (\mbox{sym} \;\kappa_g)_{2 \times 2}\right)~\mbox{d}x',
\end{split}
\end{equation}
whose two integral terms are strictly tied to the deformation of the geometric mid-surface of $\Omega^h$, with the first term representing stretching and the second term the bending both relative to the growth tensor(see Remark 1 at the end of section \ref{VariableT_FvK1} for more heuristics).

In section \ref{VariableT_FvK2}, we compute the Euler-Lagrange equations associated with $\mathcal I_g$ in (\ref{Ig_intro}), in the case of isotropic materials. These equations are expressed in terms of the Airy stress potential $\Phi$, Young's modulus $S$, Poisson's ratio $\nu$ and the bending stiffness $B$:
\begin{equation}\label{EL_intro}
\begin{split}
& \frac{1}{g_1 + g_2} \Delta^2 \Phi + \zeta(\Phi) = -S\left(K_G +\lambda_g\right)\\
& B(g_1 + g_2)^3 \Delta^2 v = (g_1 + g_2)[\Phi, v] + (\nabla(g_1 + g_2))^T\mbox{cof} \nabla^2 \Phi \nabla v - B \Omega_g - B \eta(v) + \frac{1}{2}\xi(\Phi).
\end{split}
\end{equation}
Also, $\lambda_g$ and $\Omega_g$ are similar to those introduced in \cite{Lewicka-Mahadevan-Pakzad_2011} modified by the thickness functions, while $\zeta(\Phi)$, $\eta(v)$ and $\xi(\Phi)$ are new terms unique to the varying thickness case, see section \ref{VariableT_FvK2} .

Finally, in section \ref{VariableT_FvK3} we establish convergence of equilibria, i.e. convergence of critical points of $I^h$ to critical points of  $\mathcal I_g$, under certain additional assumptions (\ref{m-w-ass}). When the material is isotropic, the set of solutions to (\ref{EL_intro}) coincides with the set of the critical points of $\mathcal I_g$.

\section{The Gamma-Convergence}\label{VariableT_FvK1}
In this section, we study the asymptotic behaviour of a deformations sequence whose energy scales of order $h^4$.
 Recall the definition of $s^h$ in (\ref{ChangeVarFvK}). Then we have:
\begin{theorem}\label{CompactnessFvK}
Assume the energies of a sequence of deformations $u^h \in W^{1, 2}(\Omega^h, \mathbb R^3)$ satisfy:
\begin{equation}\label{EnergyBoundFvK}
I^h(u^h) \leq C h^4,
\end{equation}
with some constant $C>0$. Then there exist rotations $\bar R^h \in 
SO(3)$ and translations $c^h \in \mathbb R^3$ such that for the normalized 
deformations:
\begin{equation}\label{NVyFvK}
y^h(x', x_3) = (\bar R^h)^T u^h\big(x', s^h(x', x_3)\big) - c^h: 
\Omega^{\ast} \to\mathbb R^3,
\end{equation}
the following assertions hold:
\begin{itemize}
\item [(i)] $y^h(x', x_3)$ converge in $W^{1, 2}(\Omega^{\ast}, \mathbb R^3)$ to 
$x'$.
\item [(ii)] The rescaled average displacements:
\begin{equation}\label{SDVFvK}
V^h(x') = \frac{1}{h} \fint_{-1/2}^{1/2} y^h(x', t) - \big(x', s^h(x', t)\big)^T
 \mathrm{d}t 
\end{equation}
converge (up to a subsequence) in $W^{1, 2}(\Omega, \mathbb R^3)$ to the vector 
field of the form $(0, 0, v)^T$, with the only non-zero out-of-plane scalar 
component: $v \in W^{2, 2}(\Omega, \mathbb R)$. 
\item [(iii)] The scaled in-plane displacements $h^{-1}V_{tan}^h$ converge  weakly in $W^{1, 2}(\Omega, \mathbb R^2)$, up to a subsequence, to an in-plane 
displacement field $w \in W^{1, 2}(\Omega, \mathbb R^2)$. 
\item [(iv)] The scaled energies $\frac{1}{h^4}I^h(u^h)$ satisfy the lower bound:
\begin{equation*}
\liminf_{h \to 0} \frac{1}{h^4} I^h(u^h) \geq \mathcal I_g(w, v),
\end{equation*}
where:
\begin{equation}\label{LMEVFvK}
\begin{split}
\mathcal I_g(w, v) =& \frac{1}{2}\int_{\Omega}(g_1 + g_2) \mathcal Q_2\Big(\mathrm{sym} \nabla w + \frac{1}{2} \nabla v \otimes \nabla v - (\mathrm{sym} \;\epsilon_g)
_{2 \times 2}\\
&\hspace{1.5cm}-\frac{1}{2}(g_2-g_1)(\mathrm{sym}\; \kappa_g)_{2 \times 2} + \frac{1}{2}\mathrm{sym}(\nabla v \otimes \nabla(g_2 - g_1))_{2 \times 2}\Big)~\mathrm{d}x'\\
&+ \frac{1}{24}\int_{\Omega}(g_1 + g_2)^3\mathcal Q_2\left(\nabla^2 v + (\mathrm{sym}\; \kappa_g)_{2 \times 2}\right)\mathrm{d}x',
\end{split}
\end{equation}
and the quadratic nondegenerate form $\mathcal Q_2$, acting on matrices $F \in 
\mathbb R^{2 \times 2}$ is:
\begin{equation}\label{VQ2Q3FvK}
\mathcal Q_2(F) = \min\{\mathcal Q_3(\tilde F)\mid \tilde F \in \mathbb R^{3 \times 3}, \tilde F_{2 \times 2} = F\}\mbox{ where }\mathcal Q_3(\tilde F) = \nabla^2W(\mathrm{Id}_3)(\tilde F, \tilde F).
\end{equation}
\end{itemize}
\end{theorem}
We anticipate that, in addition to the compactness analysis above, we further prove existence of a recovery sequence which realizes the lower bound in (iv), namely:
\begin{theorem}\label{RecoveryFvK}
For every $w \in W^{1, 2}(\Omega, \mathbb R^2)$ and every $v \in W^{2, 2}(\Omega, 
\mathbb R)$, there exists a sequence of deformations $u^h \in W^{1, 2}(\Omega^h, 
\mathbb R^3)$ such that the following holds as $h\to 0$:
\begin{itemize}
\item [(i)] The sequence $y^h(x', x_3) = u^h(x', s^h(x', x_3))$ converges in 
$W^{1, 2}(\Omega^{\ast}, \mathbb R^3)$ to $x'$.
\item [(ii)] $V^h(x')$ defined as in (\ref{SDVFvK}) converge in $W^{1, 2}(\Omega, 
\mathbb R^3)$ to $(0, 0, v)^T$.
\item [(iii)] $h^{-1}V_{tan}^h$ converge in $W^{1, 2}(\Omega, \mathbb R^2)$ to $w$.
\item [(iv)] The limit of the corresponding scaled energies  
realizes (\ref{LMEVFvK}):
\[
\lim_{h \to 0}\frac{1}{h^4}I^h(u^h) = \mathcal I_g(w, v).
\]  
\end{itemize}
\end{theorem}
An essential ingredient in the proof of Theorem \ref{CompactnessFvK} is the following approximation lemma, obtained through the geometric rigidity estimate  in \cite{Friesecke-James-Muller_2002}:
\begin{lemma}\label{GeoAppVPFvK}
Let $u^h \in W^{1, 2}(\Omega^h, \mathbb R^3)$ satisfy:
\[
\lim_{h \to 0}\frac{1}{h^2}I^h(u^h) = 0.
\]
Then there exist matrix fields $R^h \in W^{1, 2}(\Omega, \mathbb R^{3 \times 3})$, 
such that $R^h(x') \in SO(3)$ for a.e. $x' \in \Omega$ and:
\begin{equation}\label{Esti1FvK}
\frac{1}{h}\int_{\Omega^h}\left|\nabla u^h(x)- R^h(x')a^h(x)\right|^2~\mathrm{d}x \leq
 C\big(I^h(u^h) + h^4\big),
\end{equation}
\begin{equation}\label{Esti2FvK}
\int_{\Omega}|\nabla R^h|^2 ~\mathrm{d}x' \leq Ch^{-2}\big(I^h(u^h) + h^4\big),
\end{equation}
with constant $C$ independent of $h$.
\end{lemma}
The proof is exactly the same as the proof of Theorem 1.6 in \cite{Lewicka-Mahadevan-Pakzad_2011}, where the Friesecke, James and M\"uller's inequality is applied to small cylinders in $\Omega^h$, hence we omit it.
Owing to Lemma \ref{GeoAppVPFvK}, there follows the compactness and lower bound part or Theorem \ref{CompactnessFvK} :
\begin{proof} 

[Theorem \ref{CompactnessFvK}]
{\bf 1.} By (\ref{EnergyBoundFvK}), (\ref{Esti1FvK})(\ref{Esti2FvK}), for each $u^h$ there exists $R^h \in W^{1, 2}(\Omega, SO(3))$ with:
\begin{equation}\label{3.2'FvK}
\frac{1}{h}\int_{\Omega^h}\big|\nabla u^h - R^h a^h\big|^2 \leq Ch^4, \quad\quad
\int_{\Omega}\big|\nabla R^h\big|^2 \leq Ch^2.
\end{equation}
Define the averaged rotations by projecting onto $SO(3)$:
\begin{equation*}
\tilde R^h = \mathbb P_{SO(3)}\fint_{\Omega}R^h.
\end{equation*}
Based on the estimate on $\nabla R^h$ in (\ref{3.2'FvK}), these projections are well defined for small $h$. Moreover:
\begin{equation}\label{3.3'FvK}
\int_{\Omega}|R^h - \tilde R^h|^2 \leq C \Big(\int_{\Omega}\big|R^h - \fint_\Omega R^h\big|^2 + \mbox{dist}^2\Big(\fint_\Omega R^h, SO(3)\Big)\Big)\leq C \int_{\Omega}|\nabla R^h|^2 \leq C h^2.
\end{equation}
Now, a further projection:
\begin{equation}\label{3.4'FvK}
\hat R^h = \mathbb P_{SO(3)}\fint_{\Omega^h}(\tilde R^h)^T \nabla u^h
\end{equation}
is also well defined for small $h$, since $\mbox{dist}^2\Big(\fint_{\Omega^h}(\tilde R^h)^T \nabla u^h, SO(3)\Big)$ is bounded by:
\begin{equation}\label{3.5'FvK}
\begin{split}
& \left|\fint_{\Omega^h} (\tilde R^h)^T \nabla u^h~\mbox{d}x' - \mbox{Id}_3\right| \leq C \fint_{\Omega^h}|\nabla u^h - \tilde R^h|^2~\mbox{d}x\\
& \qquad\qquad \leq C\left(\fint_{\Omega^h}|\nabla u^h - R^h a^h|^2 + \fint_{\Omega^h}|a^h - \mbox{Id}_3|^2 + \fint_{\Omega^h}|R^h - \tilde R^h|^2\right) 
\leq C h^2.
\end{split}
\end{equation}
Consequently, we obtain:
\begin{equation}\label{3.6'FvK}
|\hat R^h - \mbox{Id}_3|^2 \leq C\Big|\fint_{\Omega^h}(\tilde R^h)^T\nabla u^h ~\mbox{d}x - \mbox{Id}_3\Big|^2 \leq C h^2.
\end{equation}

{\bf 2.} Define a new approximating rotation in:
\begin{equation}\label{3.7'FvK}
\bar R^h = \tilde R^h \hat R^h.
\end{equation}
This will be the final rotation in the definition (\ref{NVyFvK}). According to (\ref{3.2'FvK}), (\ref{3.3'FvK}) and (\ref{3.6'FvK}):
\begin{equation}\label{3.8'FvK}
\int_{\Omega}|R^h - \bar R^h|^2 \leq C h^2 \quad \mbox{and} \quad \lim_{h \to 0}(\bar R^h)^T R^h = \mbox{Id} \quad \mbox{in}~~ W^{1, 2}(\Omega, \mathbb R^{3 \times 3}).
\end{equation}
Choose $c^h \in \mathbb R^3$ such that for the rescaled average displacement $V^h$ in (\ref{SDVFvK}), we have:
\begin{equation}\label{3.91FvK}
\int_{\Omega} V^h = 0.
\end{equation}
Since for any $F$ sufficiently close to $SO(3)$, its projection $\dps{ R = \mathbb P_{SO(3)} F}$ coincides with the rotation appearing in the polar decomposition $F = RU$ where $\mbox{skew}\, U = 0$, it follows that $U = R^T F = (\mathbb P_{SO(3)} F)^T F$ is symmetric. In the virtue of (\ref{3.4'FvK}) and (\ref{3.7'FvK}), this implies that:
\[
(\bar R^h)^T \fint_{\Omega^h} \nabla u^h = (\hat R^h)^T (\tilde R^h)^T \fint_{\Omega^h} \nabla u^h = \left(\mathbb P_{SO(3)} \fint_{\Omega^h} (\tilde R^h)^T \nabla u^h \right)^T (\tilde R^h)^T \fint_{\Omega^h} \nabla u^h
\] 
is symmetric as well. On the other hand, $\fint_{\Omega^h} \nabla u^h$ is close to $\tilde R^h$, hence to $SO(3)$, in virtue of (\ref{3.5'FvK}).
Together with the above equality, this observation implies:
\begin{equation*}
\bar R^h = \mathbb P_{SO(3)} \fint_{\Omega^h} \nabla u^h.
\end{equation*}
We will next calculate the gradient of the normalized deformation $y^h$, then apply Poincar\'e's inequality to prove (i). From (\ref{3.2'FvK}) and (\ref{3.5'FvK}), we get:
\begin{equation*}
\begin{split}
\|\nabla_{x'}&y^h - (\mbox{Id}_3)_{3 \times 2}\|_{L^2(\Omega^{\ast})}^2 \\
&\leq \int_{\Omega^{\ast}} \left|(\bar R^h)^T \Big(\nabla_{tan} u^h(x', s^h) + \partial_3 u^h(x', s^h(x', x_3))\otimes \nabla_{x'} s^h\Big) - (\mbox{Id}_3)_{3 \times 2}\right|^2 ~\mbox{d}x\\
& \leq C \frac{1}{h}\int_{\Omega^h} |(\bar R^h)^T \nabla u^h - \mbox{Id}_3|^2 + Ch \leq Ch,
\end{split}
\end{equation*}
and also:
\begin{equation*}
\begin{split}
\|\partial_3 y^h\|_{L^2(\Omega^{\ast})}^2 & =  \int_{\Omega} \int_{-1/2}^{1/2}(g_1^h + g_2^h)^2 \left|(\bar R^h)^T \partial_3 u^h(x', s^h(x', x_3))\right|^2 ~\mbox{d}x_3 \mbox{d}x'\\
&\leq Ch \int_{\Omega^h} |\partial_3 u^h|^2 \leq Ch \int_{\Omega^h}|\nabla u^h|^2 \leq Ch.
\end{split} 
\end{equation*}
In conclusion:
\begin{equation}\label{3.102FvK}
\nabla y^h \to \nabla x' \quad\quad \mbox{in} ~ L^2(\Omega^{\ast}).
\end{equation}
Observe that the choice of $c^h$ gives us that:
\[
0 = \int_{\Omega} V^h = \frac{1}{h}\int_{\Omega^{\ast}}\big(y^h(x', x_3) - x'\big) -\int_{\Omega^{\ast}}\big[0, \frac{1}{h}s^h(x', x_3)\big]^T, 
\]
which further implies:
\[
\Big|\int_{\Omega^{\ast}} y^h(x', x_3)~\mbox{d}x - x' \Big| = h \left|\int_{\Omega^{\ast}}\big[0, \frac{1}{h}s^h(x', x_3)\big]\right| \leq Ch \to 0.
\]
Application of Poincar\'e's inequality finally yields (i), because:
\[
\begin{split}
\|y^h - x' \|_{L^2(\Omega^{\ast})} & \leq \big\|y^h - x' - \fint_{\Omega^{\ast}}(y^h - x')\big\|_{L^2(\Omega^{\ast})} + \big\|\fint_{\Omega^{\ast}}y^h - x'\big\|_{L^2(\Omega^{\ast})}\\
& \leq C\|\nabla y^h - \nabla x'\|_{L^2(\Omega^{\ast})} + Ch \to 0.
\end{split}
\]

{\bf 3.} Consider the matrix fields $A^h \in W^{1,2}(\Omega, \mathbb R^{3 \times 3})$ defined as:
\begin{equation}\label{3.11'FvK}
\begin{split}
A^h(x') & = \frac{1}{h}\fint_{-g_1^h}^{g_2^h} \left((\bar R^h)^T R^h(x')a^h(x', t) - \mbox{Id}_3\right)~\mbox{d}t  \\ & = \frac{1}{h}(\bar R^h)^T R^h(x') \Big(\fint_{-g_1^h}^{g_2^h} a^h(x', t) ~\mbox{d} t \Big) \mbox{Id}_3\\
& = \frac{1}{h}\left((\bar R^h)^T R^h(x') - \mbox{Id}_3\right) + h(\bar R^h)^T R^h(x') \epsilon_g(x') + \frac{1}{2}(g_2^h - g_1^h)(\bar R^h)^T R^h(x') \kappa_g.
\end{split}
\end{equation}
Thanks to (\ref{3.8'FvK}), (\ref{3.2'FvK}) and to the properties of $g_1^h, g_2^h$, we get that ${\|A^h\|_{W^{1. 2}(\Omega)}\leq C}$, and so:
\begin{equation}\label{3.12'FvK}
\begin{split}
\lim_{h \to 0} A^h &= A \quad \mbox{and} \quad \lim_{h \to 0}\frac{1}{h}\left((\bar R^h)^T R^h - \mbox{Id}_3\right) = A, \\ 
& \mbox{weakly in }~W^{1,2}(\Omega, \mathbb R^{3 \times 3}) \quad \mbox{and (strongly) in }~ L^q(\Omega, \mathbb R^{3 \times 3}) \quad \forall q \geq 1,
\end{split}
\end{equation}
up to a subsequence. Again, applying (\ref{3.8'FvK}) and (\ref{3.2'FvK}), results in:
\[
\begin{split}
& \frac{1}{h}\left\|\mbox{sym}\big((\bar R^h)^T R^h - \mbox{Id}_3\big)\right\|_{L^2(\Omega)} = \frac{1}{2h}\left\|\big((\bar R^h)^T R^h - \mbox{Id}_3\big)^T\big((\bar R^h)^T R^h - \mbox{Id}_3\big)\right\|_{L^2(\Omega)} \\
& \qquad\qquad\qquad \leq C \frac{1}{h} \|(\bar R^h)^T R^h - \mbox{Id}_3 \|_{L^4(\Omega)}^2 \leq  C \frac{1}{h} \|R^h - \bar R^h \|_{W^{1, 2}(\Omega)}^2 \leq Ch. 
\end{split}
\]
Thus, the limiting matrix field $A$ is skew-symmetric:
\begin{equation}\label{3.14'FvK}
\mbox{sym}\, A = \lim_{h \to 0} \mbox{sym} \frac{1}{h}
((\bar R^h)^T R^h - \mbox{Id}_3) = 0.
\end{equation} 
In addition, we notice that:
\[
\begin{split}
&\frac{1}{h} \mbox{sym}\, A^h = \mbox{sym}\big((\bar R^h)^T R^h \epsilon_g(x')\big) + \frac{1}{2h}({g_2^h} - {g_1^h})\,\mbox{sym}\left((\bar R^h)^TR^h \kappa_g\right) \\
&\hspace{6cm}- \frac{1}{2h^2}\left((\bar R^h)^TR^h - \mbox{Id}_3\right)^T \left((\bar R^h)^TR^h - \mbox{Id}_3\right)
\end{split} 
\]
Therefore, the properties of $g_1^h, g_2^h$, (\ref{3.8'FvK}), (\ref{3.12'FvK}) and (\ref{3.14'FvK}) imply:
\begin{equation}\label{3.15'FvK}
\lim_{h \to 0} \frac{1}{h} \mbox{sym}\,A^h = \mbox{sym}\,\epsilon_g + \frac{1}{2}(g_2 - g_1)\mbox{sym} \,\kappa_g + \frac{1}{2} A^2 \quad \mbox{in}~ L^q(\Omega, \mathbb R^{3 \times 3}) \quad \forall q \geq 1. 
\end{equation}

{\bf 4.} Concerning the convergence of $V^h$, a direct calculation indicates:
\begin{equation}\label{3.16'FvK}
\begin{split}
\nabla V^h(x') = A^h_{3 \times 2}(x') & + \frac{1}{h}(\bar R^h)^T \int_{-1/2}^{1/2}\left(\nabla_{tan}u^h(x', s^h(x', t))- R^h(x')a^h(x', s^h(x', t)\right))~\mbox{d}t 
\\ & + \frac{1}{h}\int_{-1/2}^{1/2}\left((\bar R^h)^T \partial_3 u^h(x', s^h(x', t)) - e_3\right) \otimes \nabla_{x'} s^h(x', t)~\mbox{d}t.
\end{split}
\end{equation}
From (\ref{3.2'FvK}), the second term in the right hand side above is bounded by $Ch$ in $L^2(\Omega)$. We can rewrite the third term in the right hand side of (\ref{3.16'FvK}) as:
\begin{equation}\label{3.161FvK}
\begin{split}
\frac{1}{h}&\int_{-1/2}^{1/2}\left((\bar R^h)^T \partial_3 u^h(x', s^h(x', t)) - e_3\right) \otimes \nabla_{x'} s^h(x', t)~\mbox{d}t\\
&= \frac{1}{h}\int_{-1/2}^{1/2}(\bar R^h)^T\Big(\nabla u^h(x', s^h(x', t)) - R^h a^h(x', s^h(x', t)) \\ & \qquad\qquad \qquad \qquad \quad + R^h\left(a^h(x', s^h(x', t)) - \mbox{Id}_3\right)
+ R^h - \bar R^h\Big)e_3 \otimes \nabla_{x'}s^h(x', t)~\mbox{d}t,
\end{split}
\end{equation}
Based on the convergence properties of of $g_1^h, g_2^h$, the definition of $a^h$, and (\ref{3.2'FvK}) and (\ref{3.8'FvK}), the third term in (\ref{3.16'FvK}) or the quantity in (\ref{3.161FvK}) is also bounded by $Ch$ in $L^2(\Omega)$. Hence we have:
\begin{equation}\label{3.162FvK}
\|\nabla V^h - A^h_{3 \times 2}\|_{L^2(\Omega)} \leq Ch.
\end{equation}
By (\ref{3.12'FvK}), the matrix field $\nabla V^h$ thus converges in $L^2(\Omega, \mathbb R^{3 \times 2})$ to $A_{3 \times 2}$. In view of (\ref{3.91FvK}), this convergence, together with Poincar\'e's inequality, implies:
\begin{equation}\label{3.17'FvK}
\lim_{h \to 0} V^h = V \quad \mbox{in} ~ W^{1, 2}(\Omega, \mathbb R^3) \qquad \mbox{and}\qquad \nabla V = A_{3 \times 2}.
\end{equation}
Since $A \in W^{1, 2}(\Omega, \mathbb R^{3 \times 3})$, we see that there must be $V \in W^{2, 2}(\Omega, \mathbb R^3)$. But $\mbox{sym} \nabla (V_{tan}) = 0$ according to  (\ref{3.14'FvK}), whereas Korn's inequality yields $V_{tan}$ being constant, thus $0$ in virtue of (\ref{3.91FvK}). This concludes the proof of (ii). For (iii), we apply Poincar\'e's and Korn's inequalities, in:
\begin{equation}\label{3.171FvK}
\begin{split}
& \left\|h^{-1}V_{tan}^h\right\|_{W^{1, 2}(\Omega)} \leq C \big\|\nabla(h^{-1} V_{tan}^h)\big\|_{L^2(\Omega)}\\
& \qquad\qquad \leq C \Big\|\nabla (h^{-1}V_{tan}^h) - h^{-1} \fint_{\Omega} \mbox{skew} \nabla V_{tan}^h \Big\|_{L^2(\Omega)}+ C \Big\|h^{-1} \fint_{\Omega} \mbox{skew} \nabla V_{tan}^h\Big\|_{L^2(\Omega)}\\
& \qquad\qquad \leq C \big\|\mbox{sym}\nabla (h^{-1}V_{tan}^h)\big\|_{L^2(\Omega)} +  C \big\|h^{-1} \fint_{\Omega} \mbox{skew} \nabla V_{tan}^h\big\|_{L^2(\Omega)} \leq C,
\end{split}
\end{equation}
where we also utilized (\ref{3.15'FvK}), (\ref{3.162FvK}), (\ref{3.12'FvK}) and the estimate for the last two terms of (\ref{3.16'FvK}). This indeed yields (iii).

{\bf 5.} Define the scaled strains $Z^h \in L^2(\Omega^{\ast}, \mathbb R^{3 \times 3})$ by setting:
\[
Z^h(x', x_3) = \frac{1}{h^2}\left((R^h(x'))^T \nabla u^h(x', s^h(x', x_3))a^h(x', 
s^h(x', x_3))^{-1}-\mbox{Id}_3\right).
\]
Owing to (\ref{3.2'FvK}), these are bounded: $\|Z^h\|_{L^2(\Omega^{\ast})} \leq C$ and hence, up to a subsequence:
\begin{equation}\label{VGhFvK}
\lim_{h \to 0} Z^h = Z \quad \quad \mbox{weakly in } \; L^2(\Omega^{\ast}, \mathbb R^{3 \times 3}).
\end{equation}
Properties for the limiting strain $Z$ are derived as follows. First observe that:
\begin{equation}\label{yh3aFvK}
\lim_{h \to 0} \frac{1}{h^2}\Big(\partial_3 y^h - (g_1^h + g_2^h) e_3\Big) = (g_1+ g_2)A e_3 \quad \quad \mbox{in} ~ L^2(\Omega^{\ast}, \mathbb R^3). 
\end{equation}
One may refer to \cite{Lewicka-Mahadevan-Pakzad_2011} for the detailed calculation, which is the same here. Second, for each small $s > 0$ we define the family of functions $f^{s, h} \in W^{1, 2}(\Omega^{\ast}, \mathbb R^3)$ in:
\begin{equation}\label{fshFvK}
\begin{split}
f^{s, h} (x) & = \frac{1}{h^2} \frac{1}{s} \left(y^h(x + se_3) - y^h(x) - (g_1^h + g_2^h)se_3\right)\\ 
& = \frac{1}{h^2}\fint_0^s \partial_3 y^h(x + t e_3) - (g_1^h + g_2^h)e_3 ~\mbox{d}t.
\end{split}
\end{equation}
By (\ref{yh3aFvK}) there holds:
\begin{equation}\label{lfshFvK}
\lim_{h \to 0} f^{s, h} = (g_1 + g_2)Ae_3\quad \mbox{ and } \quad \lim_{h \to 0} \partial_3 f^{s, h} = 0 \quad \quad \mbox{in}~ L^2(\Omega^{\ast}, \mathbb R^3),
\end{equation}
because:
\[
\partial_3 f^{s, h}(x)= \frac{1}{s} \frac{1}{h^2}\left(\partial_3 y^h(x + se_3) - \partial_3 y^h(x)\right).
\]
Further, for any $\alpha = 1, 2$, we have:
\[
\begin{split}
\partial_{\alpha}&f^{s, h}(x) = \frac{1}{h^2}\frac{1}{s} \Bigg((\bar R^h)^T\Big(\partial_{\alpha} u^h(x', s^h(x', x_3+s))- \partial_{\alpha} u^h(x', s^h(x', x_3))\Big)\\
&+(\bar R^h)^T\Big(\partial_{\alpha}s^h(x', x_3 + s)\partial_3 u^h(x', s^h(x', x_3+s))- \partial_{\alpha}s^h(x', x_3)\partial_3 u^h(x', s^h(x', x_3))\Big)\\
&\hspace{10cm}\qquad  - \partial_{\alpha}(g_1^h + g_2^h)se_3\Bigg).
\end{split}
\]
We split the right hand side above into two parts and investigate them separately. The first term:
\begin{equation}\label{3.211FvK}
\begin{split}
& \frac{1}{h^2}\frac{1}{s}(\bar R^h)^T\Big(\partial_{\alpha} u^h(x', s^h(x', x_3+s))- \partial_{\alpha} u^h(x', s^h(x', x_3))\Big)\\
& = (\bar R^h)^T R^h\Bigg(\frac{1}{s}\left(Z^h(x', x_3 + s) - Z^h(x', x_3)\right)a^h\left(x', s^h(x', x_3 + s)\right)\\
&\hspace{7cm} + \left(\mbox{Id}_3 + h Z^h(x', x_3)\right)\frac{g_1^h + g_2^h}{h}\Bigg)e_{\alpha}
\end{split}
\end{equation}
weakly converges in $L^2(\Omega^{\ast})$ to:
\[
\left(\frac{1}{s}\left(Z(x', x_3+s) - Z(x', x_3)\right) + (g_1 + g_2)\kappa_g\right)e_{\alpha},
\]
by (\ref{VGhFvK}) and the properties of $g_1^h, g_2^h$. The second part can be rewritten as:
\[
\begin{split}
& \frac{1}{h^2}\frac{1}{s}\Bigg((\bar R^h)^T\left(\partial_{\alpha}s^h(x', x_3 + s)\partial_3 u^h (x', s^h(x', x_3 + s)) - \partial_{\alpha} s^h(x', x_3)\partial_3 u^h(x', s^h(x', x_3))\right)\\
&\hspace{10.5cm} - s\partial_{\alpha}(g_1^h + g_2^h)e_3\Bigg)\\
&= \frac{1}{h^2}\frac{1}{s}\Big((\bar R^h)^T \partial_3 u^h(x', s^h(x', x_3 + s))\big(\partial_{\alpha}s^h(x', x_3 + s) - \partial_{\alpha} s^h(x', x_3)\big) - s\partial_{\alpha}(g_1^h + g_2^h)e_3\Big)\\
&\quad  + \frac{1}{h^2}\frac{1}{s} \partial_{\alpha}s^h(x', x_3)(\bar R^h)^T\left(\partial_3 u^h(x', s^h(x', x_3 + s)) - \partial_3 u^h(x', s^h(x', x_3))\right).
\end{split}
\]
Using the previously derived estimates (\ref{3.2'FvK}) and (\ref{3.8'FvK}), we obtain:
\[
\begin{split}
&\frac{1}{h^2}\frac{1}{s}\Big((\bar R^h)^T \partial_3 u^h(x', s^h(x', x_3 + s))\left(\partial_{\alpha}s^h(x', x_3 + s) - \partial_{\alpha} s^h(x', x_3)\right) - s\partial_{\alpha}(g_1^h + g_2^h)e_3\Big)\\
&=\frac{1}{h^2}\partial_{\alpha}(g_1^h + g_2^h)(\bar R^h)^T\left(\nabla u^h(x', s^h(x', x_3 + s)) - \bar R^h\right)e_3\\
&=\frac{1}{h^2}\partial_{\alpha}(g_1^h + g_2^h)(\bar R^h)^T\Bigg(\nabla u^h(x', s^h(x', x_3 + s)) - R^h a^h(x', s^h(x', x_3 + s)) \\
&\hspace{5cm}+ R^h \left(a^h(x', s^h(x', x_3 + s)) - \mbox{Id}_3\right) + R^h - \bar R^h\Bigg)e_3\\
&\to \partial_{\alpha}(g_1 + g_2)Ae_3 \quad \mbox{in} ~ L^2(\Omega^{\ast}),
\end{split}
\]
and further:
\[
 \frac{1}{h^2}\frac{1}{s} \partial_{\alpha}s^h(x', x_3)(\bar R^h)^T\left(\partial_3 u^h(x', s^h(x', x_3 + s)) - \partial_3 u^h(x', s^h(x', x_3))\right) \to 0 \quad \mbox{in}~ L^2(\Omega^{\ast}).
\]
Hence, in view of the above analysis, there follows:
\begin{equation}\label{3.21'FvK}
\lim_{h \to 0}\partial_{\alpha}f^{s, h}(x) = \frac{1}{s}\left(Z(x', x_3 + s) - Z(x', x_3)\right)e_{\alpha} + (g_1 + g_2) \kappa_g e_{\alpha} + \partial_{\alpha}(g_1 + g_2)Ae_3,
\end{equation}
weakly in $L^2(\Omega^{\ast})$. Consequently, $f^{s, h}$ converges weakly in $W^{1, 2}(\Omega^{\ast}, \mathbb R^3)$ to $(g_1 + g_2)Ae_3$. Equating the tangential derivatives, we thus obtain:
\[
\partial_{\alpha}\left((g_1 + g_2)Ae_3\right) = \frac{1}{s}\left(Z(x', x_3 + s) - Z(x', x_3)\right)e_{\alpha} + (g_1 + g_2) \kappa_g e_{\alpha} + \partial_{\alpha}(g_1 + g_2)Ae_3, 
\]
for $\alpha = 1, 2$, which is equivalent to:
\begin{equation}\label{3.22'FvK}
Z(x', x_3)e_{\alpha} = Z(x', 0)e_{\alpha} + x_3(g_1 + g_2)Z_1(x')e_{\alpha},
\end{equation}
where:
\begin{equation}\label{3.23'FvK}
Z_1(x') = \nabla(Ae_3) - \kappa_g.
\end{equation}

{\bf 6.} We will now calculate $\mbox{sym}Z(x', 0)_{2 \times 2}$, through computing $1/h \,\mbox{sym} \nabla V^h$. Divide both sides of (\ref{3.16'FvK}) by $h$ and observe that the second term there can be rewritten as:
\[
\begin{split}
&\frac{1}{h^2}(\bar R^h)^T\int_{-1/2}^{1/2}\left(\nabla_{tan}u^h(x', s^h(x', t)) - R^h(x')a^h(x', s^h(x', t))\right)~\mbox{d}t\\
& = \frac{1}{h^2}(\bar R^h)^TR^h \int_{-1/2}^{1/2}\left((R^h)^T \nabla_{tan}u^h(x', s^h(x', t))a^h(x', s^h(x', t))^{-1} - \mbox{Id}_3\right)a^h(x', s^h(x', t))~\mbox{d}t\\
&= (\bar R^h)^TR^h \int_{-1/2}^{1/2}Z^h(x', t)a^h(x', s^h(x', t))~\mbox{d}t\\ 
&= (\bar R^h)^TR^h \int_{-1/2}^{1/2}Z^h(x', t)\left(\mbox{Id}_3 + h^2 \epsilon_g + hs^h(x', t)\kappa_g\right)\mbox{d}t. 
\end{split}
\]
Thus, weakly in $L^2(\Omega)$, there exists the following limit:
\[
\begin{split}
&\lim_{h \to 0} \mbox{sym} \frac{1}{h^2}(\bar R^h)^T\int_{-1/2}^{1/2}\left(\nabla_{tan}u^h(x', s^h(x', t)) - R^h(x')a^h(x', s^h(x', t))\right)~\mbox{d}t \\
& = \lim_{h \to 0} \mbox{sym} \left((\bar R^h)^TR^h \int_{-1/2}^{1/2}Z^h(x', t)\left(\mbox{Id}_3 + h^2 \epsilon_g + hs^h(x', t)\kappa_g\right)~\mbox{d}t\right)= \mbox{sym}Z(x', 0).
\end{split}
\]
Divide both sides of (\ref{3.16'FvK}) by $h$ and pass to the weak limit with the symmetric parts:
\begin{equation}\label{3.231FvK}
\begin{split}
\lim_{h \to 0} \frac{1}{h} \mbox{sym} \nabla V^h = & \; \mbox{sym} \, \epsilon_g + \frac{1}{2}(g_2 - g_1)\mbox{sym}\,\kappa_g + \frac{1}{2}A^2 + \mbox{sym}Z(x', 0) \\
& - \frac{1}{2} \mbox{sym}\big(\nabla v \otimes (\nabla_{x'}(g_2 - g_1))\big).
\end{split}
\end{equation}
Meanwhile, by (iii), $1/h \, \mbox{sym} \nabla V^h_{tan} \to \mbox{sym}\nabla w$ weakly in $L^2(\Omega, \mathbb R^{2 \times 2})$. Consequently:
\begin{equation}\label{G0FvK}
\begin{split}
\mbox{sym}\,Z(x', 0)_{2 \times 2}
= & \; \mbox{sym} \nabla w - (\mbox{sym} \,\epsilon_g)_{2 \times 2} - \frac{1}{2}(g_2 - g_1)(\mbox{sym} \,\kappa_g)_{2 \times 2} - \frac{1}{2}(A^2)_{2 \times 2}\\
& + \frac{1}{2}\mbox{sym}\big( \nabla v \otimes (\nabla_{x'}(g_2 - g_1))\big).
\end{split}
\end{equation}

{\bf 7.} In this final step we shall prove the lower bound in (iv). By change of variables we get:
\begin{equation}\label{IhCVFvK}
I^h(u^h) = \frac{1}{h} \int_{\Omega^h} W(\nabla u^h (a^h)^{-1}) = \int_{\Omega}\frac{g_1^h + g_2^h}{h}\int_{-1/2}^{1/2}W\big(\mbox{Id}_3+h^2Z^h(x', x_3)\big)~\mbox{d}x_3\mbox{d}x'.
\end{equation}
On the "good" set $\Omega_h = \{x \in \Omega^{\ast}\mid h|Z^h(x', x_3)| \leq 1\}$ we use the Taylor expansion:
\begin{equation}\label{FormalTaylorFvK}
\frac{1}{h^4}W\big(\mbox{Id}_3 + h^2 Z^h(x', x_3)\big) = \frac{1}{2} \mathcal Q_3\big(Z^h(x', x_3)\big) + o(1)|Z^h|^2.
\end{equation} 
On the other hand, the characteristic functions $\chi_h$ of $\Omega_h$ satisfy:
\begin{equation}\label{chihFvK}
\lim_{h \to 0}\chi_h = 1 \qquad \mbox{in}~ L^1(\Omega^{\ast}),
\end{equation}
as $hZ^h$ converges to $0$ pointwise a.e. by (\ref{3.2'FvK}). Hence, there follows:
\[
\begin{split}
\liminf_{h \to 0} \frac{1}{h^4} I^h(u^h) & \; \geq \liminf_{h \to 0} \int_{\Omega^{\ast}}\chi_h\frac{g_1^h + g_2^h}{h} W\big(\mbox{Id}_3 + h^2 Z^h(x', x_3)\big)~\mbox{d}x\\
& = \liminf_{h \to 0} \int_{\Omega^{\ast}}\chi_h \frac{g_1^h + g_2^h}{h}\left(\frac{1}{2}\mathcal Q_3\big((Z^h(x', x_3)\big) + o(1)|Z^h|^2\right)~\mbox{d}x\\
& \geq \liminf_{h \to 0}\frac{1}{2}\int_{\Omega^{\ast}}\frac{g_1^h + g_2^h}{h}\mathcal Q_3(\chi_h Z^h) ~\mbox{d}x = \frac{1}{2} \int_{\Omega^{\ast}}(g_1 + g_2)\mathcal Q_3\left(Z(x', x_3)\right)\mbox{d}x_3\mbox{d}x'.
\end{split}
\]
Since $\mathcal Q_3$ is positive definite and depends only on the symmetric part of its argument, we get:
\[
\begin{split}
\mathcal Q_3(Z(x', x_3)) & \; = \mathcal Q_3 \big(\mbox{sym}\,Z(x', x_3)\big)\geq \mathcal Q_2\big(\mbox{sym}\,Z(x', x_3)_{2 \times 2}\big)\\
&= \mathcal Q_2\big(\mbox{sym}\,Z(x', 0)_{2 \times 2} + x_3 (g_1 + g_2) Z_1(x')\big) \\
& = \mathcal Q_2\left(\mbox{sym} \,Z(x', 0)_{2 \times 2}\right) + (g_1 + g_2)^2 x_3^2 \mathcal Q_2\left(\mbox{sym}\,Z_1(x')_{2 \times 2}\right)\\
& \quad + 2 x_3 (g_1 + g_2)\mathcal L_2\big(\mbox{sym}\,Z(x', 0)_{2 \times 2}, \mbox{sym}\, Z_1(x')_{2 \times 2}\big),
\end{split}
\]
where $\mathcal L_2$ is the corresponding bilinear form of $\mathcal Q_2$. Therefore:
\[
\begin{split}
&\liminf_{h \to 0} \frac{1}{h^4}I^h(u^h)\\
&\geq \frac{1}{2} \int_{\Omega}(g_1 + g_2)\int_{-1/2}^{1/2}\mathcal Q_2\big(\mbox{sym} \;Z(x', 0)_{2 \times 2}\big) + (g_1 + g_2)^2 x_3^2 \mathcal Q_2\big(\mbox{sym}\; Z_1(x')_{2 \times 2}\big)~\mbox{d}x_3\mbox{d}x'\\
& = \frac{1}{2}\int_{\Omega}(g_1 + g_2)\mathcal Q_2\big(\mbox{sym} \; Z(x', 0)_{2 \times 2}\big)~\mbox{d}x' +\frac{1}{24} \int_{\Omega}(g_1+g_2)^3 \mathcal Q_2\big(\mbox{sym}\; Z_1(x')_{2 \times 2}\big)~\mbox{d}x',
\end{split}
\]
which implies that:
\[
\begin{split}
&\liminf_{h \to 0} \frac{1}{h^4}I^h(u^h)\\
& \geq \frac{1}{2}\int_{\Omega}(g_1 + g_2)\mathcal Q_2\Big(\mbox{sym} \nabla w - (\mbox{sym}\;\epsilon_g)_{2 \times 2} - \frac{1}{2}(A^2)_{2 \times 2} - \frac{1}{2}(g_2 - g_1)(\mbox{sym}\; \kappa_g)_{2 \times 2}\\
& \qquad\qquad\qquad \qquad \qquad +\frac{1}{2}\mbox{sym}\left(\nabla v \otimes \nabla(g_2 - g_1)\right)\Big)~\mbox{d}x'\\ 
&  \quad + \frac{1}{24}\int_{\Omega}(g_1 + g_2)^3\mathcal Q_2\Big(\mbox{sym}(\nabla(Ae_3) - \kappa_g)_{2 \times 2}\Big)~\mbox{d}x'.
\end{split}
\]
In view of (ii) and (\ref{3.17'FvK}), we note that:
\[
(A^2)_{2 \times 2} = - \nabla v \otimes \nabla v \qquad \mbox{and} \qquad Ae_3 = - \nabla v.
\]
This concludes the proof of (iv) and of the Theorem.
\end{proof}

In the remaining part of this section, we will present the crucial points of proving Theorem \ref{RecoveryFvK}. For more detailed proof, we refer to \cite{Lewicka-Mahadevan-Pakzad_2011}. To construct a recovery sequence with claimed properties, for any $F \in \mathbb R^{2 \times 2}$, let $(F)^{\ast} \in \mathbb R^{3 \times 3}$ denote the matrix for which $(F)^{\ast}_{2 \times 2} = F$ and $(F)_{i3}^{\ast} = (F)_{3i}^{\ast} = 0$ for $ i = 1, 2, 3$. Also, let $c(F)\in \mathbb R^3$ be the unique vector satisfying $\mathcal Q_2(F) = \mathcal Q_3\left((F)^{\ast} + \mbox{sym}(c \otimes e_3)\right)$. The well-definedness and the linearity of the mapping $c: \mathbb R_{\mathrm{sym}}^{2 \times 2} \to \mathbb R^3$ is due to the positive definiteness of the quadratic form $\mathcal Q_3$ on the space of symmetric matrices. We also need to set $l(F)$ for all $F \in \mathbb R^{3 \times 3}$ to be the unique vector in $\mathbb R^3$, such that:
\[
\mbox{sym}\big(F - (F_{2 \times 2})^{\ast}\big) = \mbox{sym}\big(l(F)\otimes e_3\big).
\]
Now, for any in-plane displacement $w$ and out-of-plane displacement $v$ as in Theorem \ref{RecoveryFvK}, their corresponding recovery sequence is given by:
\begin{equation}\label{VRecFvK}
u^h(x', x_3) = \begin{bmatrix}x'\\ x_3\end{bmatrix} + \begin{bmatrix}h^2 w(x') \\ h v(x')\end{bmatrix} + \Big(x_3 - \frac{1}{2}(g_2^h - g_1^h)\Big)\begin{bmatrix}-h \nabla v(x') \\ 0\end{bmatrix} + h^2 x_3 d^0(x') + \frac{1}{2} hx_3^2 d^1(x'), 
\end{equation}
where:
\begin{equation}\label{vd1d2FvK}
\begin{split}
& d^0 = l(\epsilon_g) + c\Big(\mbox{sym} \nabla w -(\mbox{sym} \; \epsilon_g)_{2 \times 2} + \frac{1}{2}\nabla v \otimes \nabla v - \frac{1}{2}(g_2 - g_1)(\mbox{sym}\; \kappa_g)_{2 \times 2}\\
& \qquad\qquad \qquad \quad + \frac{1}{2}\mbox{sym}(\nabla v \otimes \nabla (g_2 - g_1))\Bigg) - \frac{1}{2}(g_2 - g_1)c\Big(- \nabla^2 v -(\mbox{sym} \;\kappa_g)_{2\times2}\Big)\\
& d^1 = l(\kappa_g) + c\Big(-\nabla^2 v - (\mbox{sym} \;\kappa_g)_{2 \times 2}\Big).
\end{split}
\end{equation}
This ends the sketch of the proof. $\hfill \Box$

\bigskip

\noindent{\bf Remark.}
The two terms in the limiting  energy $\mathcal I_g(w, v)$ in (\ref{LMEVFvK}) are strictly tied to the deformations of the geometric mid-surface $\tilde{\phi}^h(\Omega)$ of $\Omega^h$. Namely, define $\tilde{\phi}^h: \Omega \to \mathbb R^3$ as: 
\[
\tilde{\phi}^h(x') = \begin{bmatrix}x'\\ \frac{1}{2}(g_2^h(x')-g_1^h(x')) \end{bmatrix},
\]
and consider the deformation: 
\[
\phi_1^h(x') = \tilde{\phi}^h(x') + \begin{bmatrix}
    h^2 w(x')\\ h v(x').
\end{bmatrix}
\]
We have:
\[
\nabla \tilde{\phi}^h = \begin{bmatrix}
    1 &0\\
    0 &1\\
    \frac{1}{2} \partial_1(g_2^h-g_1^h) &\frac{1}{2} \partial_1(g_2^h-g_1^h)
\end{bmatrix},\quad 
\nabla \phi_1^h = \begin{bmatrix}
    1+h^2 \partial_1 w_1 &h^2 \partial_2 w_1\\
    h^2 \partial_1 w_2 &1+h^2 \partial_2 w_2\\
    \frac{1}{2} \partial_1(g_2^h-g_1^h)+h\partial_1 v &\frac{1}{2} \partial_2(g_2^h-g_1^h)+h\partial_2 v
\end{bmatrix}.
\]
Given $\tau \in T_x(\Omega)$, the change of the first fundamental form of $\tilde{\phi}^h(\Omega)$ equals:
\[
\begin{split}
&\left|\partial_{\tau} \phi_1^h \right|^2 - \left|\left(a^h\circ \tilde{\phi}^h\right) \left(\partial_{\tau} \tilde{\phi}^h\right)\right|^2\\  
& = 2h^2 \tau^T \Big(\mathrm{sym} \nabla w + \frac{1}{2} \nabla v \otimes \nabla v - (\mathrm{sym} \; \epsilon_g)_{2 \times 2}
- \frac{1}{2}\left(g_2-g_1\right)(\mathrm{sym} \; \kappa_g)_{2 \times 2}\\
&\qquad\qquad \qquad +\frac{1}{2}\mathrm{sym}\left(\nabla v \otimes \nabla \left(g_2 - g_1\right)\right)\Big) \tau + \mathcal O (h^3).
\end{split}
\]
Hence, the expression in the argument of $\mathcal Q_2$ in the first term of (\ref{LMEVFvK}) describes {\it stretching}, namely the second order in $h$ change of the first fundamental form of the geometric mid-surface $\tilde{\phi}^h(\Omega)$ in relation to the growth tensor $a^h$.  

To understand the second term of $\mathcal I_g(w, v)$, we consider the change of the second fundamental form of  $\tilde{\phi}^h(\Omega)$ in relation to $a^h$. For each $\tau, \eta \in T_x(\Omega)$, we want to estimate the difference:
\begin{equation}\label{Hbending}
    \left \langle \Pi^h \partial_{\tau} \phi_1^h, \partial_{\eta} \phi_1^h \right \rangle - \Big \langle \big(\frac{1}{2} \partial_3 G^h + \tilde{\Pi}^h\big)\partial_{\tau} \tilde{\phi}^h, \partial_{\eta} \tilde{\phi}^h\Big \rangle,
\end{equation}
where $\Pi^h$ is the shape operator on $\phi_1^h(\Omega)$, and $\tilde \Pi^h$ is the shape operator on $\tilde{\phi}^h(\Omega)$, and where $G^h = (a^h)^T a^h$ is the Riemannian metric corresponding to the growth tensor $a^h$. The first term in (\ref{Hbending}) measures the bending of the deformed geometric mid-surface $\phi^h_1(\Omega)$. The second term measures the bending of the geometric mid-surface plus the bending effect of the Riemannian metric induced by the growth tensor $a^h$. To better understand the bending effect of $G^h$, we refer to Remark 11.8 (ii) on page 279 of \cite{Lewicka_2023}.

Similar to the analysis in Remark 4.3 of \cite{Lewicka-Mora-Pakzad_2010}, we then have:
\[
\tilde{\Pi}^h \partial_{\tau} \tilde{\phi}^h = \tilde{\Pi}^h(\nabla \tilde{\phi}^h)\tau = \partial_{\tau}\left(\frac{\tilde{\bf n}^h}{\left|\tilde{\bf n}^h\right|}\right), \qquad \Pi^h \partial_{\tau} \phi_1^h =\Pi^h(\nabla \phi_1^h)\tau = \partial_{\tau}\left(\frac{{\bf n}_1^h}{\left|{\bf n}_1^h\right|}\right),
\]
where $\tilde{\bf n}^h = \partial_1 \tilde{\phi}^h \times \partial_2 \tilde{\phi}^h$ is the unit normal of $\tilde{\phi}^h(\Omega)$, and ${\bf n}_1^h = \partial_1 \phi_1^h \times \partial_2 \phi_1^h$ is unit normal of $\phi_1^h(\Omega)$. Through straightforward calculation, we obtain:
\[
\tilde{\bf n}^h = \begin{bmatrix}
    -\frac{1}{2} \partial_1 (g_2^h - g_1^h)\\
    -\frac{1}{2} \partial_2 (g_2^h - g_1^h)\\
    1
\end{bmatrix}, \qquad
{\bf n}_1^h = \begin{bmatrix}
    -\frac{1}{2} \partial_1 (g_2^h - g_1^h) - h \partial_1 v\\
    -\frac{1}{2} \partial_2 (g_2^h - g_1^h) - h \partial_2 v\\ 
    1
\end{bmatrix} + \mathcal{O}(h^2),
\]
so that, in particular, $|\tilde{\bf n}^h| = 1+ \mathcal{O}(h^2)$ and $|{\bf n}_1^h|= 1+ \mathcal{O}(h^2)$. Therefore:
\[
\begin{split}
&\tilde{\Pi}^h(\nabla \tilde{\phi}^h)\tau = \partial_{\tau} \tilde{\bf n}^h + \mathcal{O}(h^2) = -\frac{1}{2} \begin{bmatrix}
    \nabla^2(g_2^h - g_1^h)\\ 0
\end{bmatrix}\tau + \mathcal{O}(h^2),
\\ &
\Pi^h\left(\nabla \phi_1^h\right)\tau = \partial_{\tau} {\bf n}_1^h + \mathcal{O}(h^2) = \begin{bmatrix}
    -\frac{1}{2} \nabla^2(g_2^h - g_1^h) - h \nabla^2 v \\
    0
\end{bmatrix}\tau + \mathcal{O}(h^2).
\end{split}
\]
Recall that:
\[
\frac{1}{2} \partial_3 G^h = h \mathrm{sym} \;\kappa_g + \mathcal{O}(h^3).
\]
The above implies that (\ref{Hbending}) equals:
\[
\begin{split}
\left \langle\Big((\nabla \phi_1^h)^T \Pi^h (\nabla \phi_1^h) - (\nabla \tilde{\phi}^h)^T\big(\frac{1}{2}\partial_3 G^h - \tilde{\Pi}^h\big) (\nabla \tilde{\phi}^h ) \Big) \tau, \eta \right\rangle \\
= - h \left \langle \left(\nabla^2 v + (\mathrm{sym}\; \kappa_g)_{\mathrm{tan}}\right)\tau, \eta \right \rangle + \mathcal{O}(h^2).
\end{split}
\]
We see that the second term of (\ref{LMEVFvK}) relates to {\it bending}, specifically the first order in $h$ change in the second fundamental form of the geometric mid-surface in relation to the growth tensor $a^h$. $\hfill\Box$

\section{The F\"oppl-von K\'arm\'an Equations}\label{VariableT_FvK2}

In this section, we will derive the Euler-Lagrange equations of the limiting energy $\mathcal{I}_g$ as in (\ref{LMEVFvK}) in case of variable thickness isotropic plates, namely under the additional property of:
\begin{equation}\label{isotropicFvK}
\forall F \in \mathbb R^{3 \times 3} \quad \forall R \in SO(3) \quad\quad W(FR) = W(F).
\end{equation}
For each $F \in \mathbb{R}^{3 \times 3}$ and $\tilde{F} \in \mathbb{R}^{2 \times 2}$, the quadratic forms $\mathcal Q_3, \mathcal Q_2$ have the expression (see e.g. \cite{Friesecke-James-Muller_2006}):
\begin{equation}\label{isoQ3FvK}
\mathcal Q_3(F) = 2 \mu|\mbox{sym}\,F|^2 + \lambda |\mbox{Tr}\, F|^2,\qquad \mathcal Q_2(\tilde F) = 2\mu|\mbox{sym} \, \tilde F|^2 + \frac{2\mu\lambda}{2\mu + \lambda}|\mbox{Tr} \, \tilde F|^2,
\end{equation}
where $\mu$ and $\lambda$ are the Lam\'e constants. 
Following the same calculation as in \cite{Lewicka-Mahadevan-Pakzad_2011}, we obtain the following Euler-Lagrange equations for (\ref{LMEVFvK}):
\begin{equation}\label{ELequationFvK}
\begin{split}
& \displaystyle{\frac{1}{(g_1 + g_2)} \Delta^2 \Phi + \zeta(\Phi) = -S\left(K_G +\lambda_g\right)},\\
& \displaystyle{B(g_1 + g_2)^3 \Delta^2 v = (g_1 + g_2)[\Phi, v] + (\nabla(g_1 + g_2))^T\mbox{cof} \,\nabla^2 \Phi \nabla v - B \Omega_g - B \eta(v) + \frac{1}{2}\xi(\Phi)}.
\end{split}
\end{equation}
We now explain the quantities above:\vspace{1mm}
\begin{itemize}
\item $~\dps{S = -\frac{\mu(2\mu + 3 \lambda)}{\mu + \lambda}}$ is Young's modulus,
\item $~\dps{\nu = \frac{\lambda}{2(\lambda + \mu)}}$ is Poisson's ratio,\vspace{1mm}
\item $~\dps{B = \frac{S}{12(1-\nu^2)}}$ is bending stiffness,\vspace{1mm}
\item $~\dps{K_G = \frac{1}{2}[v, v] = \det \nabla^2 v}$ is the Gaussian curvature of the deformed midsurface,
\item $~\dps{\lambda_g = \mbox{curl}^T\mbox{curl} \Big((\epsilon_g)_{2 \times 2} - \frac{1}{2}(g_2-g_1)(\mbox{sym} \;\kappa_g)_{2 \times 2} + \frac{1}{2} \nabla v \otimes \nabla(g_2 - g_1)\Big)}$.\vspace{1mm}
\item $~\dps{\zeta(\Phi) =  2\nabla\Big(\frac{1}{g_1 + g_2}\Big)\cdot\nabla(\Delta \Phi) + \frac{S}{2\mu}\nabla^2\Big(\frac{1}{g_1+g_2}\Big):\nabla^2 \Phi - \nu\Delta \Big(\frac{1}{g_1 + g_2}\Big)\Delta \Phi}$,\vspace{1mm}
\item $~\dps{\eta(v) = \big(\nabla ((g_1 + g_2)^3)\big)^T\mbox{div}\nabla^2 v +  \nabla^2\big((g_1 + g_2)^3)\big):(\nabla^2 v + \nu \,\mbox{cof}\, \nabla^2 v)}$,\vspace{1mm}
\item $~\dps{\xi(\Phi) = (g_1 + g_2)[\Phi, g_2 - g_1] + (\nabla(g_2 + g_1))^T(\mbox{cof} \, \nabla^2 \Phi)\nabla(g_2 - g_1)}$,\vspace{1mm}
\item $~\dps{\Omega_g = \big\langle\nabla^2(g_1+g_2)^3:\big((\mbox{sym} \;\kappa_g)_{2 \times 2} + \nu \,\mbox{cof}\,(\mbox{sym} \;\kappa_g)_{2 \times 2}\big)\big\rangle + \nabla (g_1 + g_2)^3 \cdot \mbox{div}((\mbox{sym} \;\kappa_g)_{2 \times 2})}$.
\end{itemize}
The Airy stress potential $\Phi \in W^{2, 2}(\Omega, \mathbb R)$ plays as a medium for recovering $w$:
\[
\begin{split}
&\mbox{cof} \,\nabla^2 \Phi = (g_1 + g_2)\left(2\mu\big(\mbox{sym} \nabla w + \Psi(v)\big) + \frac{2 \mu \lambda}{2\mu + \lambda}\big(\mbox{div} w + \mbox{Tr}\,\Psi(v)\big)\mbox{Id}_2\right)
\\ &
\mbox{where}\qquad \Psi(v) = \frac{1}{2} \nabla v \otimes \nabla v - (\mbox{sym} \;\epsilon_g)_{2 \times 2} + \frac{1}{2} \mbox{sym}\big(\nabla v \otimes \nabla(g_2 - g_1)\big) - \frac{1}{2}(\mbox{sym} \; \kappa_g)_{2 \times 2},
\end{split}
\]
and the Airy bracket $[\cdot, \cdot]$ is defined as: $[v, \Phi] = \big\langle\nabla^2 v : (\mbox{cof} \nabla^2 \Phi)\big\rangle$. 

\smallskip

The natural boundary conditions associated with (\ref{ELequationFvK}) are:
\begin{equation}\label{BdryFvK}
\begin{split}
&\Phi = \partial_{\vec n} \Phi = 0,\\
&\langle \tilde \Psi : (\vec n \otimes \vec n)\rangle + \langle\nu \tilde \Psi :(\tau \otimes \tau)\rangle = 0, \hspace{7cm} \mbox{on}~\partial \Omega,\\
& (1-\nu) \partial_{\tau}\Big\langle(g_1 + g_2)^3 \tilde \Psi : (\vec n \otimes \tau)\big\rangle + \mbox{div}\Big((g_1 + g_2)^3 (\tilde \Psi + \nu \, \mbox{cof} \,\tilde \Psi)\Big)\vec n = 0.
\end{split}
\end{equation}
where $\tilde \Psi = \nabla^2 v + (\mbox{sym} \,\kappa_g)_{2 \times 2}$, and where $\vec{n}, \tau$ denote the unit normal and the unit tangent to $\partial \Omega$, respectively.
In particular, when $g_1 = g_2 = 1/2$, the system (\ref{ELequationFvK}), (\ref{BdryFvK}) coincides with the one obtained in \cite{Lewicka-Mahadevan-Pakzad_2011}.

\section{Convergence of Equilibria}\label{VariableT_FvK3}

In this section, we consider the convergence of equilibria under physical growth conditions for the energy density. As in \cite{Mora-Scardia_2012}, the density $W: \mathbb R^{3 \times 3} \to [0, +\infty]$, in addition to (\ref{w-ass}), shall satisfy:
\begin{equation}\label{m-w-ass}
\left\{
\begin{minipage}{12.5cm}
\begin{itemize}
\item [(v)] $W$ is of class $\mathcal C^1$ on $\mathbb R_{+}^{3 \times 3}$ of $3 \times 3$ matrices with positive determinant.
\item [(vi)] $W(F) = +\infty$ if $\det F \leq 0$, and $W(F) \to +\infty$ as $\det F \to 0+$.
\item [(vii)] $|\nabla W(F)F^T|\leq C(W(F) + 1)$ for every $F \in \mathbb R_+^{3 \times 3}$ and some uniform $C>0$.  
\end{itemize} 
\end{minipage}
\right.
\end{equation}
Here, the growth requirement  (vii) is assumed for $\nabla W$, and it is compatible with the blow-up requirement (vi), as pointed out in \cite{Ball_2002}. Besides, due to (vi), one cannot legitimately perform the external variation $u^h + \varepsilon \phi$ of a minimizer $u^h$ to obtain the Euler-Lagrange equations in the conventional weak form \cite{Ball_2002}. Instead, we shall consider the internal variations $u^h + \varepsilon \phi \circ u^h$, whereas the equilibrium condition for $u^h$ becomes:
\begin{equation}\label{StaFvK}
\int_{\Omega^h} \Big\langle\nabla W\Big(\nabla u^h (a^h)^{-1}\Big)\Big(\nabla u^h (a^h)^{-1}\Big)^T : \nabla \phi(u^h)\Big\rangle~ \mbox{d}x = 0 \qquad \forall \phi \in \mathcal C_b^1(\mathbb R^3, \mathbb R^3).
\end{equation}
We refer to $u^h$ as the stationary point of the energy $I^h$, if (\ref{StaFvK}) is satisfied. The space $\mathcal C_b^1$ consists of the bounded $\mathcal C^1$ functions.
We will also use the bilinear form $\mathcal{L}_2$ associated with $\mathcal Q_2$, which has already been used in the proof of Theorem \ref{CompactnessFvK}. More precisely:
\[
\mathcal L_2(E, F) = \frac{1}{2}\Big(\mathcal Q_2(E + F) - \mathcal Q_2(E) - \mathcal Q_2(F)\Big) \qquad \forall E, F \in \mathbb R^{2 \times 2}.
\]
Since $\mathcal Q_2$ depends only on the symmetric part of its argument, we have:
\begin{equation}\label{LsymFvK}
\mathcal L_2(E, F) = \mathcal L_2(\mbox{sym}\;E, \mbox{sym}\;F) = \mathcal L_2(\mbox{sym}\;E, F) = \mathcal L_2(E, \mbox{sym}\;F).
\end{equation}
With a small abuse of notation, for each $E \in \mathbb R^{2 \times 2}$ we define a linear functional $\mathcal L_2 E$ on $\mathbb R^{2 \times 2}$, by setting: $\langle \mathcal L_2 E:F \rangle = \mathcal L_2(E, F)$ for each $F \in \mathbb R^{2 \times 2}$.

Calculating the variations of $\mathcal I_g(w, v)$ in $w$ and $v$ respectively, we obtain the following weak formulation of the Euler-Lagrange equations for $\mathcal{I}_g$ as in (\ref{LMEVFvK}):
\begin{equation}\label{EL1FvK}
\begin{split}
&\int_{\Omega} (g_1+g_2)\Big\langle\mathcal L_2\big(\mbox{sym} ~\nabla w + \frac{1}{2}
\nabla v \otimes \nabla v - (\mbox{sym}\;
\epsilon_g)_{2 \times 2} -\frac{1}{2}(g_2 - g_1)\left(\mbox{sym}\; \kappa_g\right)_{2 \times 2} \\
&\hspace{5cm}+ \frac{1}{2} \mbox{sym}\big(\nabla \otimes \nabla(g_2 - g_1)\big)_{2 \times 2}\big):\mbox{sym} \nabla \psi \Big\rangle ~\mbox{d}x' = 0,
\end{split}
\end{equation}
and:
\begin{equation}\label{EL2FvK}
\begin{split}
&\int_{\Omega} (g_1 + g_2)\Big\langle\mathcal L_2\big(\mbox{sym} \; \nabla w + \frac{1}{2}
\nabla v \otimes \nabla v - (\mbox{sym}~
\epsilon_g)_{2 \times 2} -\frac{1}{2}(g_2 - g_1)\left(\mbox{sym}\; \kappa_g\right)_{2 \times 2} \\
&\hspace{3cm}+ \frac{1}{2} \mbox{sym}\big(\nabla \otimes \nabla(g_2 - g_1)\big)_{2 \times 2}\big):\big(\nabla v + \frac{1}{2} \nabla(g_2 - g_1)\big) \otimes \nabla \varphi\Big\rangle ~\mbox{d}x'\\
& \quad + \frac{1}{12} \int_{\Omega} (g_1 + g_2)^3 \Big\langle\mathcal L_2\big(\nabla^2 v + (\mbox{sym}\;
\kappa_g)\big):\nabla^2 \varphi \Big\rangle ~\mbox{d}x' = 0,
\end{split}
\end{equation}
for any $\psi \in \mathcal C_b^1(\mathbb R^2, \mathbb R^2)$ and $\varphi
\in C_b^2(\mathbb R^2)$.

\bigskip

The stated convergence of equilibria is contained in the following result:
\begin{theorem}\label{MainTFvK}
Assume $u^h \in W^{1, 2}(\Omega^h; \mathbb R^3)$ to be a sequence of stationary points of $I^h$ with
\begin{equation}\label{energysclFvK}
I^h(u^h) \leq C h^4.
\end{equation}
Then there exist  $\bar R^h \in SO(3)$ and
 $c^h \in \mathbb R^3$, such that for the normalized
deformations:
$$ y^h(x', x_3) = (\bar R^h)^T u^h\big(x', s^h(x',x_3)\big) - c^h: \Omega^{\ast}\to \mathbb R^3, $$
there hold the convergence properties (i), (ii) and (iii) in Theorem \ref{CompactnessFvK} and moreover:
\begin{itemize}
\item [(iv)] $(v, w)$ solves (\ref{EL1FvK}) and (\ref{EL2FvK}).
\end{itemize}
\end{theorem}
The proof of the theorem is based on the method presented in \cite{Muller-Pakzad_2008} and developed in \cite{Lewicka_2011, Mora-Scardia_2012}. The following is our detailed proof. 

\begin{proof}
{\bf 1.} As before, (\ref{energysclFvK}) implies (i), (ii) and (iii) of Theorem \ref{CompactnessFvK},  Also, based on (\ref{3.15'FvK}):
\begin{equation}\label{estiFvK}
\left\|(\bar R^h)^T \nabla u^h(x', s^h(x', x_3)) - \mbox{Id}_3\right\|^2_{L^2(\Omega^{\ast})} \leq C \fint_{\Omega^h} \left|(\bar R^h)^T \nabla u^h(x', x_3) - \mbox{Id}_3\right|^2 ~\mbox{d}x\leq Ch^2.
\end{equation}
Noticing that:
\[
\partial_3 y_3^h = (g_1^h + g_2^h)\left((\bar R^h)^T \nabla u^h(x', s^h(x', x_3)) \right)_{33}
\]
 and applying Poincar\'e-Wirtinger's inequality with bound (\ref{estiFvK}), we obtain:
\begin{multline*}
\Big\|\frac{y_3^h}{h} - \frac{g_1^h + g_2^h}{h}x_3 -\frac{1}{2}\frac{g_2^h - g_1^h}{h} - V_3^h(x')\Big\|_{L^2(\Omega^{\ast})} \leq C\Big\|\frac{\partial_3 y^h_3}{h}-\frac{g_1^h + g_2^h}{h}\Big\|_{L^2(\Omega^{\ast})} \\
\leq C \Big\|\frac{\partial_3 y^h_3}{g_1^h + g_2^h} - 1\Big\|_{L^2(\Omega^{\ast})} \leq \left\|(\bar R^h)^T \nabla u^h(x', s^h(x', x_3)) - \mbox{Id}_3\right\|_{L^2(\Omega^{\ast})}\leq Ch.
\end{multline*}
Together with the properties of $g_1^h, g_2^h$ and (ii), the above implies:
\begin{equation}\label{y3hFvK}
\lim_{h \to 0} \frac{y_3^h}{h} = v + (g_1 + g_2)x_3 + \frac{1}{2}(g_2-g_1) \qquad \mbox{in }~ L^2(\Omega^{\ast}).
\end{equation}
As in the proof of Theorem \ref{CompactnessFvK}, define the scaled strains $Z^h \in L^2(\Omega^{\ast}, \mathbb R^{3 \times 3})$ in:
\begin{equation}\label{FStrainFvK}
Z^h(x', x_3) = \frac{1}{h^2}\left(R^h(x')^T \nabla u^h(x', s^h(x', x_3)) a^h(x', s^h(x', x_3))^{-1} - \mbox{Id}_3\right).
\end{equation}
As before, $Z^h$ weakly converges, up to a subsequence, to some $Z$ in $L^2(\Omega^{\ast}, \mathbb R^{3 \times 3})$, satisfying:
\begin{equation}\label{limitGhFvK}
\begin{split}
& Z(x', x_3) e_{\alpha} = Z(x', 0)e_{\alpha} + x_3(g_1 + g_2)(-\nabla^2 v - \kappa_g)e_{\alpha}, \quad \mbox{for}~ \alpha = 1, 2, \\
& \mbox{where  }
\mbox{sym}\,Z(x', 0) = \mbox{sym} \nabla w - (\mbox{sym} \;\epsilon_g)_{2 \times 2} - \frac{1}{2}(g_2 - g_1)(\mbox{sym} \; \kappa_g)_{2 \times 2}+\frac{1}{2}\nabla v \otimes \nabla v \\
& \qquad\qquad\qquad \qquad \qquad\quad 
- \frac{1}{2}\mbox{sym}(\nabla v \otimes \nabla(g_2 -g_1))
\end{split}
\end{equation}

{\bf 2.} Define the scaled stress $E^h: \Omega^{\ast} \to \mathbb R^{3 \times 3}$ as:
\begin{equation}\label{EhFvK}
E^h(x', x_3) = \frac{1}{h^2} \nabla W\big(\mbox{Id}_3 + h^2 Z^h(x', x_3)\big)\big(\mbox{Id}_3 + h^2 Z^h(x', x_3)\big)^T.
\end{equation}
Such $E^h(x)$ is symmetric due to the frame indifference of $W$, and it obeys the estimate:
\begin{equation}\label{BEhFvK}
|E^h| \leq C \Big(\frac{1}{h^2}W(\mbox{Id} + h^2 Z^h) + |Z^h|\Big),
\end{equation}
and for detailed proof, one may refer to that of (4.14) in \cite{Mora-Scardia_2012} and to the argument in \cite{Lewicka_2011}.

{\bf 3.} By the definition of a stationary point of $I^h$ in (\ref{StaFvK}), we get for every $\phi \in \mathcal C_b^1(\mathbb R^3, \mathbb R^3)$
\begin{equation*}
\int_{\Omega^h}\Big\langle \nabla W\Big(\nabla u^h(x)a^h(x)^{-1}\Big)\big(\nabla u^h(x) a^h(x)^{-1})\big)^T :\nabla \phi(u^h(x))\Big\rangle~ \mbox{d}x = 0.
\end{equation*}
Using Fubini's Theorem and a change of variable, we can rewrite the above as:
\begin{equation}\label{Sta2FvK}
\begin{split}
& \int_{\Omega^{\ast}}\frac{g_1^h + g_2^h}{h} \Big\langle \nabla W\Big(\nabla u^h(x', s^h(x', x_3))\big(a^h(x', s^h(x', x_3))\big)^{-1}\Big)\cdot\\
&\hspace{.5cm}\cdot\big(\nabla u^h(x', s^h(x', x_3))\big(a^h(x', s^h(x', x_3))\big)^{-1}\big)^T:\nabla \phi\big(u^h(x', s^h(x', x_3))\big) \Big\rangle~\mbox{d}x_3 \mbox{d}x' = 0.
\end{split}
\end{equation}
For each test function $\tilde \phi \in \mathcal C_b^1(\mathbb R^3, \mathbb R^3)$ and $u \in W^{1, 2}(\Omega^h,\mathbb R^3)$, define:
\[
\phi(u) = \bar R^h\tilde \phi\big((\bar R^h)^T u - c^h\big).
\]
Recalling that $u^h = \bar R^h(y^h(x', x_3) + c^h)$ and taking the derivative, we get:
\begin{equation*}
\nabla \phi(u^h) = \bar R^h \nabla \tilde \phi\big((\bar R^h)^T u^h - c^h\big)(\bar R^h)^T = \bar R^h \nabla \tilde \phi(y^h) (\bar R^h)^T.
\end{equation*}
Substituting the above into (\ref{Sta2FvK}), we obtain that for all $\tilde \phi \in \mathcal C_b^1(\mathbb R^3, \mathbb R^3)$:
\begin{equation}\label{5.8MSFvK}
\begin{split}
&\int_{\Omega^{\ast}}\frac{g_1^h + g_2^h}{h}\Big\langle (\bar R^h)^T \nabla W \big(\nabla u^h(x',s^h(x', x_3))a^h(x', s^h(x', x_3))^{-1}\big)\cdot\\ 
&\hspace{1.5cm}\cdot \Big(\nabla u^h(x',s^h(x', x_3))a^h(x', s^h(x', x_3))^{-1}\Big)^T \bar R^h : \nabla \tilde \phi(y^h(x', x_3)) \Big\rangle ~\mbox{d}x_3 \mbox{d}x' = 0,
\end{split}
\end{equation}
Furthermore, by definition of $Z^h$ and $E^h$ in (\ref{FStrainFvK}), (\ref{EhFvK}) and by the frame indifference of $W$:
\[
\begin{split}
\nabla W&\big(\nabla u^h(x',s^h(x', x_3))a^h(x', s^h(x', x_3))^{-1}\big) \Big(\nabla u^h(x',s^h(x', x_3))a^h(x', s^h(x', x_3))^{-1}\Big)^T\\
&= R^h(x') \nabla W\big(\mbox{Id}_3 + h^2 Z^h(x', x_3)\big)\Big(\mbox{Id}_3 + h^2 Z^h(x', x_3)\Big)^T R^h(x')^T\\
& = h^2 R^h(x')E^h(x', x_3)R^h(x')^T.
\end{split}
\]
Thus, in terms of the stress $E^h$, we may rewrite (\ref{5.8MSFvK}) as:
\begin{equation}\label{SSE2FvK}
\int_{\Omega^{\ast}}\frac{g_1^h + g_2^h}{h}\Big\langle (\bar R^h)^TR^h(x')E^h(x', x_3)R^h(x')^T \bar R^h : \nabla \tilde \phi(y^h(x', x_3))\Big\rangle ~\mbox{d}x = 0.
\end{equation}

{\bf 4.} By the energy scaling (\ref{energysclFvK}), the bound (\ref{BEhFvK}) of $E^h$ and the fact that $Z^h$ are bounded in $L^2(\Omega^{\ast}, \mathbb R^{3 \times 3})$, for each measurable set $\Lambda \subset \Omega^{\ast}$, we have:
\[
\int_{\Lambda} |E^h| \mathrm{d}x 
\leq C \int_{\Lambda} \frac{1}{h^2}W(\mbox{Id}_3 + h^2 Z^h) \mathrm{d}x + C \int_{\Lambda}|Z^h|\mathrm{d}x 
\leq Ch^2 + C|\Lambda|^{1/2}.
\]
Thus, the scaled stresses $E^h$ are bounded and equi-integrable in $L^1(\Omega^{\ast}, \mathbb R^{3 \times 3})$. Hence, by the Dunford-Pettis theorem, there exists $E \in L^1(\Omega^{\ast}, \mathbb R^{3 \times 3})$ such that:
\begin{equation}\label{CEhFvK}
E^h \rightharpoonup E \quad \mbox{weakly in } ~ L^1(\Omega^{\ast}, \mathbb R^{3 \times 3}).
\end{equation}
In particular, $E$ is symmetric from the symmetry of $E^h$. In order to pass to the limit in (\ref{SSE2FvK}), a more refined convergence property of $E^h$ is necessary. Define sets:
$$B_h = \{x \in \Omega^{\ast} \mid h^{2 - \gamma}|Z^h(x)|\leq 1\},$$ 
with a chosen exponent $\gamma \in (0, 1)$ and let $\tilde \chi_h$ denote the characteristic function of $B_h$. Together with the properties of $g_1^h, g_2^h$, and following the analysis of (4.20) and (4.21) in \cite{Mora-Scardia_2012}, we obtain:
\begin{equation}\label{CEhwFvK}
\begin{split}
& (1 - \tilde \chi_h)E^h \to 0 \quad \mbox{strongly in}~L^1(\Omega^{\ast}, \mathbb R^{3 \times 3}),\\
& \tilde \chi_h E^h \rightharpoonup \mathcal L_3 Z \quad \mbox{weakly in}~L^2(\Omega^{\ast}, \mathbb R^{3 \times 3}),
\end{split}
\end{equation}
where $\mathcal L_3$ is the bilinear form corresponding to $\mathcal Q_3$. Along with the $\mathcal C_b^1$ regularity of test functions, this mixed type of convergence is sufficient for (\ref{CEhFvK}) to imply that $E = \mathcal L_3 Z \in L^2(\Omega^{\ast}, \mathbb R^{3 \times 3})$. Finally, since $(\bar R^h)^T R^h$ is bounded and converging in measure to $\mbox{Id}_3$, together with (\ref{CEhwFvK}), this yields that $\tilde \chi_h(\bar R^h)^T R^h(x')E^h(x', x_3) \rightharpoonup \mathcal L_3 Z$ weakly in $L^2(\Omega^{\ast}, \mathbb R^{3 \times 3})$.

{\bf 5.} We shall now investigate the properties of $u^h$ based on the definition of stationary points in (\ref{StaFvK}). Fix $\phi \in \mathcal C_b^1(\mathbb R^3, \mathbb R^3)$ and take $\phi^h(x) = h \phi(x', x_3/h)$, which is an admissible test function that we may insert in (\ref{SSE2FvK}), obtaining:
\[
\begin{split}
0 &= \int_{\Omega^{\ast}}\frac{g_1^h + g_2^h}{h}\Big\langle(\bar R^h)^TR^h(x')E^h(x', x_3)R^h(x')^T\bar R^h : \nabla \phi^h(y^h(x', x_3))\Big\rangle~\mbox{d}x\\
&=\int_{\Omega^{\ast}}(g_1^h + g_2^h)\sum_{\alpha = 1}^2 (\bar R^h)^TR^h(x')E^h(x', x_3)R^h(x')^T\bar R^he_{\alpha}\cdot \partial_{\alpha} \phi\big((y^h)', \frac{y_3^h}{h}\big)~\mbox{d}x \\
& \quad + \int_{\Omega^{\ast}}\frac{g_1^h + g_2^h}{h}(\bar R^h)^TR^h(x')E^h(x', x_3)R^h(x')^T\bar R^he_3 \cdot \partial_3 \phi\big((y^h)', \frac{y^h_3}{h}\big)~\mbox{d}x.
\end{split}
\]
As $(\bar R^h)^TR^h(x')E^h(x', x_3)R^h(x')^T\bar R^h$ is bounded in $L^1(\Omega^{\ast}, \mathbb R^{3 \times 3})$ and as $\partial_{\alpha} \phi$ is bounded for $\alpha =1, 2$, the first term in the right hand side of the above equality converges to zero as $h \to 0$. Thus:
\begin{equation}\label{4.24MSFvK}
\lim_{h \to 0}\int_{\Omega^{\ast}}\frac{g_1^h + g_2^h}{h}(\bar R^h)^TR^h(x')E^h(x', x_3)R^h(x')^T\bar R^he_3 \cdot \partial_3 \phi\big((y^h)', \frac{y^h_3}{h}\big)~\mbox{d}x = 0.
\end{equation}
Meanwhile, (i), (ii) and (\ref{y3hFvK}) imply:
\[
\partial_3 \phi \big((y^h)', \frac{y_3^h}{h}\big) \to \partial_3 \phi \big(x', v(x') + (g_1(x') + g_2(x'))x_3 + \frac{1}{2}(g_2(x') - g_1(x'))\big) \quad \mbox{in}~L^2(\Omega, \mathbb R^3).
\]
We now split the integral in (\ref{4.24MSFvK}) as:
\[
\begin{split}
&\int_{\Omega^{\ast}}\tilde \chi_h \frac{g_1^h + g_2^h}{h}(\bar R^h)^T R^h(x')E^h(x', x_3)R^h(x')^T \bar R^h e_3 \cdot \partial_3 \phi\big((y^h)', \frac{y_3^h}{h}\big)~\mbox{d}x\\
& + \int_{\Omega^{\ast}}(1 - \tilde \chi_h) \frac{g_1^h + g_2^h}{h}(\bar R^h)^T R^h(x')E^h(x', x_3)R^h(x')^T \bar R^h e_3 \cdot \partial_3 \phi\big((y^h)', \frac{y_3^h}{h}\big)~\mbox{d}x, 
\end{split}
\]
and apply the respective convergences of $E^h$, $g_1^h, g_2^h$ and $(\bar R^h)^T R^h$, to get:
\begin{equation}\label{4.27MSFvK}
\int_{\Omega^{\ast}}(g_1 + g_2)Ee_3\cdot \partial_3 \phi\Big(x', v + \frac{1}{2}(g_2 - g_1) + (g_1 + g_2)x_3\Big) = 0 \quad \forall \phi \in C_b^1(\mathbb R^3, \mathbb R^3).
\end{equation}
Let $v_k \in \mathcal C_b^1(\mathbb R^3)$ be a sequence of functions whose restrictions to $\Omega$ converge to $v$, strongly in $L^2(\Omega)$. Given any $\phi \in \mathcal C_b^1(\mathbb R^3, \mathbb R^3)$, we choose:
\[
\begin{split}
&\phi_k(x', x_3) = \phi\Big(x', \frac{1}{g_1 + g_2}\big(x_3 - v_k - \frac{1}{2}(g_2 - g_1)\big)\Big),\\ &
\mbox{so that }\quad 
\partial_3 \phi_k = \frac{1}{g_1 + g_2}\partial_3 \phi\Big(x', \frac{1}{g_1+ g_2}\big(x_3 - v_k - \frac{1}{2}(g_2 - g_1)\big)\Big),
\end{split}
\]
Inserting $\phi_k$ into (\ref{4.27MSFvK}), we attain:
\[
0 = \int_{\Omega^{\ast}}Ee_3 \cdot \partial_3\phi\big(x', x_3 + \frac{v - v_k}{g_1 + g_2}\big)~\mbox{d}x \to \int_{\Omega^{\ast}}Ee_3 \cdot \partial_3 \phi(x', x_3)~\mbox{d}x \quad \mbox{ as}~k \to+\infty.  
\]
Hence:
\begin{equation}
\int_{\Omega^{\ast}}Ee_3 \cdot \partial_3 \phi ~\mbox{d}x = 0 \qquad \forall \phi \in C_b^1(\mathbb R^3, \mathbb R^3).
\end{equation}
and therefore there must be $Ee_3 = 0$ a.e. in $\Omega^{\ast}$. In view of the symmetry of $E$, this implies:
\begin{equation}\label{EFvK}
E = \begin{bmatrix} E_{11} &E_{12} &0\\ E_{12} & E_{22} &0\\ 0 &0 &0 \end{bmatrix}.
\end{equation}

{\bf 6.} In this next step we investigate the zeroth moment $\bar E: S \to \mathbb R^{3 \times 3}$ of the limit stress $E$:
\begin{equation}\label{barEFvK}
\bar E(x') = \int_{-1/2}^{1/2} E(x', x_3) ~\mbox{d}x_3 \qquad \forall x' \in \Omega.
\end{equation}
We will derive the equations satisfied by $\bar E$.
To this end, consider $\psi \in \mathcal C_b^1(\mathbb R^2, \mathbb R^2)$ and choose $\tilde \phi(x) = (\psi(x'), 0)$ in (\ref{SSE2FvK}), to get:
\begin{equation}\label{4.31MSFvK}
\int_{\Omega^{\ast}}\Big\langle \frac{g_1^h + g_2^h}{h}\left[(\bar R^h)^T R^h(x')E^h(x', x_3)(R^h(x'))^h \bar R^h\right]_{2 \times 2}:\nabla \psi ((y^h)')~\mbox{d}x = 0.
\end{equation}
As in the previous step, it is convenient to split the above integral as:
\begin{equation}\label{4.32MSFvK}
\begin{split}
& \int_{\Omega^{\ast}}\tilde \chi_h\frac{g_1^h + g_2^h}{h}\Big\langle \left[(\bar R^h)^T R^h(x')E^h(x', x_3)R^h(x')^T \bar R^h\right]_{2 \times 2}:\nabla \psi ((y^h)')\Big\rangle~\mbox{d}x\\
&+ \int_{\Omega^{\ast}}(1 - \tilde \chi_h)\frac{g_1^h + g_2^h}{h}\Big\langle \left[(\bar R^h)^T R^h(x')E^h(x', x_3)R^h(x')^T \bar R^h\right]_{2 \times 2}:\nabla \psi((y^h)')\Big\rangle~\mbox{d}x.
\end{split}
\end{equation}
By (i), together with continuity and boundedness of $\nabla \psi$:
\[
\nabla \psi((y^h)') \to \nabla \psi \quad \mbox{in}~L^2(\Omega, \mathbb R^{2 \times 2}),
\]
while the weak convergence of $\tilde \chi_h E^h$, and (\ref{3.8'FvK}) imply that:
\[
\begin{split}
\lim_{h \to 0} \int_{\Omega^{\ast}} &\tilde \chi_h\frac{g_1^h + g_2^h}{h}\Big\langle \left[(\bar R^h)^T R^h(x')E^h(x', x_3)R^h(x')^T \bar R^h\right]_{2 \times 2}:\nabla \psi ((y^h)')\Big\rangle ~\mbox{d}x \\ & = \int_{\Omega^{\ast}}(g_1 + g_2) \langle E_{2 \times 2}:\nabla \psi\rangle ~ \mbox{d}x. 
\end{split}
\]
Hence, the boundedness of $\nabla \psi$ and the convergence in (\ref{CEhwFvK}) indicate that the second term in (\ref{4.32MSFvK}) converges to $0$ as $h \to 0$, and by (\ref{4.31MSFvK}), we conclude:
\[
\int_{\Omega^{\ast}}(g_1 + g_2) \langle E_{2 \times 2}:\nabla \psi \rangle ~\mbox{d}x = 0 \qquad \forall \psi \in \mathcal C_b^1(\mathbb R^2, \mathbb R^2). 
\]
The above equality can be rewritten in terms of the zeroth moment as:
\begin{equation}\label{4.33MSFvK}
\int_{\Omega}(g_1 + g_2) \langle \bar E_{2 \times 2}:\nabla \psi \rangle ~\mbox{d}x' = 0,
\end{equation}
for each $\psi \in \mathcal C_b^1(\mathbb R^2, \mathbb R^2)$, and by approximation, also for each $\psi \in W^{1, 2}(\Omega, \mathbb R^2)$.

{\bf 7.} Next, we study the equation satisfied by the first moment of stress, which is defined as:
\begin{equation}\label{hatEFvK}
\hat E(x') = \int_{-1/2}^{1/2} x_3 E(x', x_3) ~\mbox{d}x_3 \quad \forall x' \in \Omega.
\end{equation}
Let $\varphi \in \mathcal C_b^2(\mathbb R^2)$ and consider $\tilde \phi(x', x_3) = (0, \frac{1}{h} \varphi(x'))$ in (\ref{SSE2FvK}). We deduce that:
\begin{equation}\label{4.35MSFvK}
\int_{\Omega^{\ast}}\frac{1}{h}\frac{g_1^h + g_2^h}{h}\sum_{\alpha = 1}^2\left[(\bar R^h)^T R^h(x')E^h(x',x_3)R^h(x')^T \bar R^h\right]_{3\alpha}\partial_{\alpha}\varphi((y^h)')~
\mbox{d}x = 0.
\end{equation}
As in the proof of Theorem \ref{CompactnessFvK}, the matrix fields $A^h$ defined as in (\ref{3.11'FvK}) enjoy the convergence properties in (\ref{3.12'FvK}). In particular, from (ii), the limit $A$ may be written in terms of $v$ as:
\begin{equation}\label{limitAhFvK}
A = \begin{bmatrix}0 &0 &-\partial_1 v \\ 0 &0 &-\partial_2 v\\ \partial_1 v &\partial_2 v &0\end{bmatrix}.
\end{equation}
Recall that (\ref{3.11'FvK}) also implies:
\[
(\bar R^h)^T R^h(x') = \big(\mbox{Id}_3 + h A^h\big)\left(\mbox{Id}_3 + h^2 \epsilon_g(x') + \frac{1}{2}h(g_2^h - g_1^h)\kappa_g\right)^{-1} = \mbox{Id}_3 + hA^h + \mathcal O(h^2).
\]
Hence, there follows the decomposition:
\begin{equation}\label{4.36MSFvK}
\begin{split}
&\frac{1}{h}(\bar R^h)^TR^h(x')E^h(x', x_3)R^h(x')^T\bar R^h\\
&\hspace{1cm} = A^h(x')E^h(x', x_3)R^h(x')^T \bar R^h
+ E^h(x', x_3)A^h(x')^T + \frac{1}{h}E^h(x', x_3) + \mathcal O(h).
\end{split}
\end{equation}
By the bound of $E^h$ in (\ref{BEhFvK}), the convergences of $A^h$ and $\tilde \chi_h$, and the boundedness of $R^h(x')^T\bar R^h$:
\begin{equation*}
(1 - \tilde \chi_h)\left(A^h(x')E^h(x', x_3)R^h(x')^T \bar R^h + E^h(x', x_3)A^h(x')^T\right) \to 0 \quad \mbox{in}~L^1(\Omega^{\ast}, \mathbb R^{3 \times 3}),
\end{equation*}
while, by (\ref{3.12'FvK}) and by the weak convergence of $\tilde \chi_h E^h$ in $L^2(\Omega, \mathbb R^{3 \times 3})$, there follows:
\begin{equation*}
\tilde \chi_h \left(A^h(x')E^h(x', x_3)R^h(x')^T \bar R^h + E^h(x', x_3)A^h(x')^T\right) \rightharpoonup AE + EA^T,
\end{equation*}
weakly in $L^q(\Omega, \mathbb R^{3 \times 3})$ for any $q > 2$. Utilizing the last two convergences, the properties of $g_1^h, g_2^h$ and the fact that $\partial_{\alpha}\varphi((y^h)') \to \partial_{\alpha}\varphi$ in $L^p(\Omega^{\ast})$ for any $p < \infty$, we conclude that:
\begin{equation}\label{4.39MSFvK}
\begin{split}
&\int_{\Omega^{\ast}}\frac{g_1^h + g_2^h}{h}\sum_{\alpha = 1}^2\left[A^h(x')E^h(x', x_3)R^h(x')^T \bar R^h + E^h(x', x_3)A^h(x')^T\right]_{3\alpha}\partial_{\alpha}\varphi\left((y^h)'\right)~ \mbox{d}x\\
&\hspace{4cm} \to \int_{\Omega^{\ast}}(g_1 + g_2)\sum_{\alpha = 1}^{2}[AE + EA^T]_{3\alpha}\partial_{\alpha} \varphi ~\mbox{d}x \qquad \mbox{as}~h \to 0.
\end{split}
\end{equation}
Note that the expression of $A$ in (\ref{limitAhFvK}) and the structure of $E$ in (\ref{EFvK}) implies:
\[
\sum_{\alpha =1}^2[AE + EA^T]_{3 \alpha}\partial_{\alpha} \phi = \big\langle E_{2 \times 2}:(\nabla v \otimes \nabla \varphi).\big\rangle.
\]
Now, recalling the definition of $\bar E$ in (\ref{barEFvK}), there follows:
\[
\int_{\Omega^{\ast}}(g_1 + g_2)\sum_{\alpha = 1}^2[AE + EA^T]_{3\alpha} \partial_{\alpha}\varphi~\mbox{d}x = \int_{\Omega}(g_1 + g_2)\langle \bar E_{2 \times 2}:(\nabla v \otimes \nabla \varphi)\rangle ~\mbox{d}x'.
\]
Let us study (\ref{4.35MSFvK}) again. Together with (\ref{4.36MSFvK}) and (\ref{4.39MSFvK}), it clearly implies:
\begin{equation}\label{4.41MSFvK}
\begin{split}
&\lim_{h \to 0}\int_{\Omega^{\ast}}\frac{g_1^h + g_2^h}{h}\sum_{\alpha = 1}^2\left[\frac{1}{h}E^h(x', x_3)\right]_{3 \alpha}\partial_{\alpha}\varphi(y^h)')~\mbox{d}x\\
&\hspace{4cm} = - \int_{\Omega}(g_1 + g_2) \langle \bar E_{2 \times 2}:(\nabla v \otimes \nabla \varphi)\rangle ~ \mbox{d}x' \quad \forall~\varphi \in C_b^1(\mathbb R^2).
\end{split}
\end{equation}
We shall write the limit in (\ref{4.41MSFvK}) in terms of the first moment $\hat E$. The main method we use is based on the one developed in \cite{Muller-Pakzad_2008}, with a modification made in \cite{Mora-Scardia_2012}. At present, we need a new test function to take care of the varying thickness. 
Let a sequence of positive numbers $\omega_h$ satisfy:
\begin{equation}\label{4.44MSFvK}
h\omega_h \to + \infty, \qquad h^{2 - \gamma} \omega_h \to 0, \qquad \mbox{as}~ h\to 0,
\end{equation} 
where $\gamma$ is the exponent as in the definition of $B_h$. Let $\theta^h \in C_b^h(\mathbb R)$ be truncations satisfying:
\begin{equation}\label{4.45MSFvK}
\begin{split}
& \theta^h(t) = t \qquad \mbox{for}~|t|\leq \omega_h, \\ 
& |\theta^h(t)|\leq |t| \qquad \mbox{for} ~t \in \mathbb R 
\\ & \|\theta^h\|_{L^{\infty}} \leq 2 \omega_h, \qquad \left\|\frac{\mbox{d}\theta^h}{\mbox{d}t}\right\|_{L^{\infty}}\leq 2.
\end{split}
\end{equation}
For any $\eta \in C_b^1(\mathbb R^2, \mathbb R^2)$, define the admissible test function $\phi^h\in \mathcal{C}^1_b(\mathbb R^3, \mathbb R^3)$ as:
\begin{equation}\label{4.48MSFvK}
\phi^h(x) = \left(\theta^h\big(\frac{x_3}{h}\big)\eta(x'), 0\right).
\end{equation}
Substituting $\tilde \phi$ in (\ref{SSE2FvK}) by $\phi^h$  leads to:
\begin{equation}\label{4.49MSFvK}
\begin{split}
&0 = \int_{\Omega^{\ast}}\frac{g_1^h + g_2^h}{h}\theta^h\big(\frac{y_3^h}{h}\big)\Big\langle \left[(\bar R^h)^TR^h(x')E^h(x', x_3)R^h(x')^T\bar R^h\right]_{2 \times 2} : \nabla \eta ((y^h)')\Big\rangle~\mbox{d}x\\
& \quad + \int_{\Omega^{\ast}}\frac{g_1^h + g_2^h}{h}\frac{1}{h}\sum_{\alpha = 1}^2\left[(\bar R^h)^TR^h(x')E^h(x', x_3)R^h(x')^T\bar R^h\right]_{\alpha 3}\eta_{\alpha}(y^h)')\cdot \frac{\mbox{d}\theta^h}{\mbox{d}t}\big(\frac{y^h_3}{h}\big)~\mbox{d}x.
\end{split}
\end{equation}
We compute the limits of the two terms above separately. Let us begin with the first one. We study the integral in the two subdomains of the usual splitting $\Omega^{\ast} = B_h \cup (\Omega^{\ast} \setminus B_h)$. In $B_h$:
\begin{equation}\label{4.50MSFvK}
\begin{split}
\lim_{h \to 0} &\int_{\Omega^{\ast}}\tilde \chi_h \frac{g_1^h + g_2^h}{h}\theta^h\big(\frac{y_3^h}{h}\big)\Big\langle\left[(\bar R^h)^TR^h(x')E^h(x', x_3)R^h(x')^T\bar R^h\right]_{2 \times 2}:\nabla \eta ((y^h)')\Big\rangle~\mbox{d}x\\
&=\int_{\Omega^{\ast}}(g_1 + g_2)\big(v+(g_1 + g_2)x_3 + \frac{1}{2}(g_2 - g_1)\big)\langle E_{2 \times 2}:\nabla \eta(x') \rangle ~\mbox{d}x\\
& = \int_{\Omega}(g_1 + g_2)\Big\langle \Big(\big(v + \frac{1}{2}(g_2-g_1)\big)\bar E_{2 \times 2} + (g_1 + g_2)\hat E_{2 \times 2}\Big):\nabla \eta \Big\rangle~ \mbox{d}x'.
\end{split}
\end{equation}
The integral on $\Omega^{\ast}\setminus B_h$, can be estimated through (\ref{4.45MSFvK}) and the fact that the bound for $E^h$ and the definition of $B_h$ imply:
\begin{equation}\label{4.16MSFvK}
\begin{split}
\int_{\Omega^{\ast}\setminus B_h}|E^h|~\mbox{d}x & \leq C\int_{\Omega^{\ast}\setminus B_h}\frac{W(\mbox{Id}_3 + h^2 Z^h)}{h^2}~\mbox{d}x + C\int_{\Omega^{\ast}\setminus B_h}|Z^h|~\mbox{d}x\\
&\leq Ch^2 + C|\Omega^{\ast}\setminus B_h|^{1/2} \leq Ch^{2 - \gamma},
\end{split}
\end{equation}
where we also used the following inequality:
$|\Omega^{\ast}\setminus B_h| \leq h^{2(2-\gamma)}$. 
Indeed:
\[
\begin{split}
&\left|\int_{\Omega^{\ast}}(1 - \tilde \chi_h)\frac{g_1^h + g_2^h}{h}\theta^h\big(\frac{y_3^h}{h}\big)\Big\langle \left[(\bar R^h)^TR^h(x')E^h(x', x_3)R^h(x')^T\bar R^h\right]_{2 \times 2}:\nabla \eta ((y^h)'))\Big\rangle ~\mbox{d}x\right|\\
&\hspace{4cm}\leq C \omega_h\|\nabla \eta\|_{L^{\infty}}\int_{\Omega^{\ast}\setminus B_h}|E^h|\leq Ch^{2 - \gamma}\omega_h \to 0 \quad \mbox{as}~ h\to 0.
\end{split}
\]
Hence, we obtain:
\begin{equation}\label{1sttermFvK}
\begin{split}
&\lim_{h \to 0}\int_{\Omega^{\ast}}\frac{g_1^h + g_2^h}{h}\theta^h\big(\frac{y_3^h}{h}\big)\Big\langle\left[(\bar R^h)^TR^h(x')E^h(x', x_3)R^h(x')^T\bar R^h\right]_{2 \times 2}:\nabla \eta ((y^h)')\Big\rangle~\mbox{d}x\\
&\hspace{2.5cm} = \int_{\Omega}(g_1 + g_2)\Big\langle\big(v + \frac{1}{2}(g_2 - g_1)\big)\bar E_{2 \times 2} + (g_1 + g_2)\hat E_{2 \times 2} : \nabla \eta \Big\rangle~\mbox{d}x'.
\end{split}
\end{equation}
To study the second integral in (\ref{4.49MSFvK}), we split it as follows:
\begin{equation}\label{4.54MSFvK}
\begin{split}
& \int_{\Omega^{\ast}}\frac{g_1^h + g_2^h}{h^2}
\sum_{\alpha = 1}^2\left[(\bar R^h)^T R^h(x')E^h(x', x_3)R^h(x')^T\bar R^h\right]_{\alpha 3}\eta_{\alpha}((y^h)')~\mbox{d}x\\
& +\int_{\Omega^{\ast}}\frac{g_1^h + g_2^h}{h^2} \sum_{\alpha = 1}^2\left[(\bar R^h)^T R^h(x')E^h(x', x_3)R^h(x')^T\bar R^h\right]_{\alpha 3}\eta_{\alpha}((y^h)')
\cdot \Big(\frac{\mbox{d}\theta^h}{\mbox{d}t}\big(\frac{y_3^h}{h}\big)-1\Big)~\mbox{d}x.
\end{split}
\end{equation}
The second term above converges to $0$ as $h \to 0$; to prove it, we define the set $D_h = \{x \in \Omega^{\ast} \mid |y_3^h(x)| \geq h\omega_h\}$. Then, by (\ref{y3hFvK}):
\[
|D_h|\leq \omega_h^{-1}\int_{D_h}\frac{|y_3^h|}{h}~\mbox{d}x \leq \omega_h^{-1}\big\|\frac{y_3^h}{h}\big\|_{L^2(\Omega^{\ast})}\|\chi_{D_h}\|_{L^2(\Omega^{\ast})} \leq C \omega_h^{-1}|D_h|^{1/2},
\]
which implies:
\begin{equation}\label{4.61MSFvK}
|D_h| \leq C \omega_h^{-2},
\end{equation}
and we recall that applying similar method as in (\ref{4.16MSFvK}), one can get:
\begin{equation}\label{4.161MSFvK}
\int_{\Lambda}|E^h|~\mbox{d}x \leq C( h^2 + |\Lambda|^{1/2} )\qquad \forall \Lambda \subset \Omega^{\ast}.
\end{equation}
Now, the integral in the second term of (\ref{4.54MSFvK}) can be reduced to:
\[
\int_{D_h}\frac{g_1^h + g_2^h}{h^2}\sum_{\alpha = 1}^2\left[(\bar R^h)^T R^h(x')E^h(x', x_3)R^h(x')^T\bar R^h\right]_{\alpha 3}\eta_{\alpha}((y^h)')\Big(\frac{\mbox{d}\theta^h}{\mbox{d}t}\big(\frac{y_3^h}{h}\big)-1\Big)~\mbox{d}x,
\]
and owing to the properties of $g_1^h, g_2^h$ and conditions in (\ref{4.45MSFvK}), (\ref{4.61MSFvK}), (\ref{4.161MSFvK}), it is bounded by:
\[
\frac{C}{h}\Big(1 + \big\|\frac{\mbox{d}\theta^h}{\mbox{d}t}\big\|_{L^{\infty}}\Big)\|\eta\|_{L^{\infty}} \int_{D_h}|E^h|~\mbox{d}x \leq Ch + \frac{C}{h}|D_h|^{1/2}\leq Ch + \frac{C}{h \omega_h} \to 0,
\]
which proves the claimed convergence.

For the first term in (\ref{4.54MSFvK}) we observe:
\[
\begin{split}
& \lim_{h \to 0} \int_{\Omega^{\ast}}\frac{g_1^h + g_2^h}{h^2} \sum_{\alpha = 1}^2\left[(\bar R^h)^T R^h(x')E^h(x', x_3)R^h(x')^T\bar R^h\right]_{\alpha 3} \eta_i((y^h)')\cdot \frac{\mbox{d}\theta^h}{\mbox{d}t}\big(\frac{y_3^h}{h}\big)~\mbox{d}x\\
&= \lim_{h \to 0} \int_{\Omega^{\ast}}\frac{g_1^h + g_2^h}{h^2} \sum_{\alpha = 1}^2\left[(\bar R^h)^T R^h(x')E^h(x', x_3)R^h(x')^T\bar R^h\right]_{\alpha 3} \eta_\alpha((y^h)')~\mbox{d}x\\
& =  \lim_{h \to 0} \int_{\Omega^{\ast}}\frac{g_1^h + g_2^h}{h}\sum_{\alpha = 1}^2\left(A^h(x')E^h(x', x_3)R^h(x')^T\bar R^h + E^h(x', x_3)A^h(x')^T\right)\eta_{\alpha}((y^h)') ~\mbox{d}x \\
&\quad + \lim_{h \to 0} \int_{\Omega^{\ast}}\frac{g_1^h + g_2^h}{h}\sum_{\alpha = 1}^2\frac{1}{h}\left[E^h(x', x_3)\right]_{\alpha 3}\eta_{\alpha}((y^h)')~\mbox{d}x\\
& = \int_{\Omega}(g_1 + g_2)\big\langle \bar E_{2 \times 2}:(\nabla v \otimes \eta)\big\rangle~\mbox{d}x' + \lim_{h \to 0} \int_{\Omega^{\ast}}\frac{g_1^h + g_2^h}{h}\sum_{\alpha = 1}^2\frac{1}{h}E^h_{\alpha 3}(x', x_3)\eta_{\alpha}((y^h)')~\mbox{d}x.
\end{split}
\]
Substituting the above calculation and (\ref{1sttermFvK}) back into (\ref{4.49MSFvK}), we obtain:
\begin{equation}\label{4.57MSFvK}
\begin{split}
& \lim_{h \to 0}\int_{\Omega^{\ast}}\frac{g_1^h + g_2^h}{h}\sum_{\alpha = 1}^2\frac{1}{h}E^h_{\alpha 3}(x', x_3)\eta_{\alpha}((y^h)')~\mbox{d}x\\
&= -\int_{\Omega}(g_1 + g_2)\Big\langle\big(v + \frac{1}{2}(g_2 - g_1)\big)\bar E_{2 \times 2} + (g_1 + g_2)\hat E_{2 \times 2} :\nabla \eta \Big\rangle ~ \mbox{d}x'\\
& \hspace{7cm}- \int_{\Omega}(g_1 + g_2)\big\langle \bar E_{2 \times 2}:(\nabla v \otimes \eta)\big\rangle~ \mbox{d}x'.
\end{split}
\end{equation}
Applying equation (\ref{4.33MSFvK}) with $\psi = (v + 1/2(g_2-g_1))\eta$, we get:
\[
\int_{\Omega}(g_1 + g_2)\Big\langle\bar E_{2 \times 2} : \big(\nabla v + \frac{1}{2} \nabla(g_2 - g_1)\big)\otimes \eta + \big(v + \frac{1}{2}(g_2-g_1)\big)\nabla \eta\Big\rangle = 0.
\]
Using the above identity in (\ref{4.57MSFvK}), consequently yields:
\begin{equation}\label{4.58MSFvK}
\begin{split}
&\lim_{h \to 0}\int_{\Omega^{\ast}}\frac{g_1^h + g_2^h}{h}\sum_{\alpha = 1}^2\frac{1}{h}E^h_{\alpha 3}(x', x_3)\eta_{\alpha}((y^h)')~\mbox{d}x\\
&\hspace{1cm} = -\int_{\Omega}(g_1 + g_2)^2 \big\langle \hat E_{2 \times 2}:\nabla \eta \big\rangle~\mbox{d}x' + \int_{\Omega}(g_1 + g_2)\big\langle \bar E_{2 \times 2}:\frac{1}{2}\nabla(g_2 - g_1)\otimes \eta \big\rangle~\mbox{d}x'.
\end{split}
\end{equation}

{\bf 8.} Let $\varphi \in C_b^2(\mathbb R^2)$. After setting $\eta = \nabla \varphi$, we compare (\ref{4.58MSFvK}) with (\ref{4.41MSFvK}) and arrive at:
\begin{equation}\label{4.62MSFvK}
\int_{\Omega}(g_1 + g_2)\Big\langle\bar E_{2 \times 2}:\big(\nabla v + \frac{1}{2}\nabla(g_2 - g_1)\big)\otimes \nabla \varphi \Big\rangle ~\mbox{d}x' = \int_{\Omega}(g_1 + g_2)^2\big\langle \hat E_{2 \times 2}:\nabla^2 \varphi \big\rangle~\mbox{d}x'.
\end{equation}
In order to arrive at the desired equations, an explicit expression of $\bar E_{2 \times 2}$ and $\hat E_{2 \times 2}$ is necessary. Since $E = \mathcal L_3 Z$ is of the form (\ref{EFvK}), it follows that $E_{2 \times 2} = \mathcal L_2 Z_{2 \times 2}$, where we refer to Proposition 3.2 in \cite{Muller-Pakzad_2008} for details. Hence, by (\ref{limitGhFvK}), there follows:
\[
\bar E_{2 \times 2} = \mathcal L_2Z(x', 0)_{2 \times 2}, \qquad \hat E_{2 \times 2} = \frac{1}{12}(g_1 + g_2)\mathcal L_2(-\nabla^2 v - \kappa_g)_{2 \times 2}.
\]
Substituting the above into (\ref{4.33MSFvK}) and (\ref{4.62MSFvK}), in view of (\ref{limitGhFvK}) and (\ref{LsymFvK}) we conclude (iv).
\end{proof}

\bibliographystyle{siam} 
\bibliography{reference}

\end{document}